\documentclass[aap,MSNbibl,dvips]{arximspdf}
\usepackage{graphicx}
\doi{10.1214/12-AAP873} 
\volume{23}
\issue{4}
\pubyear{2013}
\firstpage{1377}
\lastpage{1408}

\makeatletter
\newcommand{\R}{\mathbb R}
\newcommand{\Z}{\mathbb Z}
\newcommand{\N}{\mathbb N}

\newtheorem{Theo}{Theorem}
\newtheorem{Lemma}[Theo]{Lemma}
\newtheorem{Cor}[Theo]{Corollary}
\newtheorem{Prop}[Theo]{Proposition}

\newproclaim{Dfn}[Theo]{Definition}

\makeatother

\begin{document}
\begin{frontmatter}

\title{Upper bound on the rate of adaptation in an~asexual population}
\runtitle{Bounding the rate of adaptation}

\begin{aug}
\author[A]{\fnms{Michael} \snm{Kelly}\corref{}\thanksref{t1}\ead[label=e1]{mbkelly@math.ucsd.edu}}
\runauthor{M. Kelly}
\affiliation{University of California, San Diego}
\address[A]{Department of Mathematics\\
University of California, San Diego\\
9500 Gilman Dr. \#0112\\
La Jolla, California 92093\\
USA\\
\printead{e1}} 
\end{aug}

\thankstext{t1}{Supported in part by NSF Grant DMS-08-05472.}

\received{\smonth{8} \syear{2011}}
\revised{\smonth{5} \syear{2012}}

%
\begin{abstract}
We consider a model of asexually reproducing individuals. The birth and
death rates of the individuals are affected by a fitness parameter. The
rate of mutations that cause the fitnesses to change is proportional to
the population size, $N$. The mutations may be either beneficial or
deleterious. In a paper by Yu, Etheridge and Cuthbertson [\textit{Ann.
Appl. Probab.} \textbf{20} (2010) 978--1004] it was shown that
the average rate at which the mean fitness increases in this model is
bounded below by $\log^{1-\delta} N$ for any $\delta> 0$. We achieve an
upper bound on the average rate at which the mean fitness increases of
$O(\log N/(\log\log N)^2)$.
\end{abstract}

%

\end{frontmatter}

\section{Introduction}\label{sec1}
In a finite, asexually reproducing population with mutations, it is
well known that competition among multiple individuals that get
beneficial mutations can slow the rate of adaptation. This phenomenon
is known as the Hill--Robertson effect, named for the authors of \cite
{HR}. One may wish to consider the effect on the rate of adaptation of
a population when there are many beneficial mutations present
simultaneously. It is easily observed that when such a population is
finite and all mutations are either neutral or deleterious, the fitness
of the population will decrease over time. This scenario is known as
Muller's ratchet. The first rigorous results regarding Muller's ratchet
were due to Haigh~\cite{H}. In an asexually reproducing population,
beneficial mutations are necessary to overcome Muller's ratchet. Yu,
Etheridge and Cuthbertson~\cite{YEC} proposed a model that gives
insight into both the Hill--Robertson effect and Muller's ratchet in
large populations with fast mutation rates.

The model introduced in~\cite{YEC} is a Moran model with mutations and
selection. There are $N$ individuals in this model, each with an
integer valued fitness. The dynamics of the model are determined by
three parameters, $\mu$, $q$ and $\gamma$, which are independent of
$N$. The parameters must satisfy $\mu> 0$, $0 < q \leq1$ and $\gamma
> 0$. Let $X_t^i$ be the fitness of individual $i$ at time $t$. Then $X
= (X^1,X^2,\ldots, X^N)$ is a stochastic process with state space
$\Z^N$. The system has the following dynamics:
\begin{longlist}[(3)]
\item[(1)] Mutation: Each individual acquires mutations at rate $\mu
$. When individual $i$ gets a mutation, it is beneficial with
probability $q$ and $X^i$ increases by 1. With probability $1-q$ the
mutation is deleterious and $X^i$ decreases by 1.
\item[(2)] Selection: For each pair of individuals $(i,j)$, at rate
$\frac{\gamma}{N}(X^i-X^j)^+$, we set $X^j$ equal to $X^i$.
\item[(3)] Resampling:\vspace*{1pt} For each pair of individuals
$(i,j)$, at rate $1/N$, we set $X^j$ equal to $X^i$.
\end{longlist}
Note that the upper bound we establish for the rate of adaptation still
holds in the absence of deleterious mutations, which corresponds to the
case $q = 1$. Under the selection mechanism the event that $X^j$ is set
to equal $X^i$ represents the more fit individual $i$ giving birth and
the less fit individual $j$ dying. Likewise, the resampling event that
causes $X^j$ to equal $X^i$ represents individual $i$ giving birth and
individual $j$ dying.

We give an equivalent description of the model involving Poisson
processes that may make the coupling arguments more clear. The Poisson
processes that determine the dynamics of $X$ are as follows:
\begin{itemize}
\item There are $N$ Poisson processes $\mathcal{P}^{i\uparrow}$, $1
\leq i \leq N$, on $[0,\infty)$ of rate $q \mu$. If $\mathcal
{P}^{i\uparrow}$ gets a mark at $t$ then the $i$th coordinate of $X$
increases by 1 at time $t$.
\item There are $N$ Poisson processes $\mathcal{P}^{i\downarrow}$, $1
\leq i \leq N$, on $[0,\infty)$ of rate $(1-q)\mu$. If $\mathcal
{P}^{i\downarrow}$ gets a mark at $t$ then the $i$th coordinate of $X$
decreases by 1 at time $t$.
\item For each ordered pair of coordinates $(i,j)$ with $i \neq j$
there is a Poisson process on $[0,\infty)$, $\mathcal{P}^{i,j}$, of
rate $1/N$. If $\mathcal{P}^{i,j}$ gets a mark at $t$ then the $j$th
coordinate changes its value to agree with the $i$th coordinate at time $t$.
\item For each ordered pair of coordinates $(i,j)$ with $i \neq j$
there is a Poisson processes on $[0,\infty) \times[0,\infty)$,
$\mathcal{P}^{i,j\uparrow}$, which has intensity $\frac{\gamma
}{N}\lambda$ where $\lambda$ is Lebesgue measure on $\R^2$. If
$\mathcal{P}^{i,j\uparrow}$ gets a mark in $\{t\} \times
[0,X_{t-}^i-X_{t-}^j]$ then the $j$th coordinate changes its value to
agree with the $i$th coordinate at time $t$.
\end{itemize}

A heuristic argument in~\cite{YEC} shows that as $N$ tends to infinity
the mean rate of increase of the average fitness of the individuals in
$X$ is $O(\log N/\break(\log\log N)^2)$. Due to a typo on page 989 they
state that the rate is $O(\log N/\break\log\log N)$. By equation (10) in
\cite{YEC},
\[
K\log(\gamma K) = 2\log N.
\]
This implies that
\[
K \approx\frac{2\log N}{\log\log N}.
\]
Plugging $2\log N/\log\log N$ into each side of the consistency
condition that they derive gives a rate of adaption of $O(\log N/(\log
\log N)^2)$.

The heuristic argument is difficult to extend to a rigorous argument. Let
\[
\overline{X} = \frac{1}{N} \sum_{i=1}^N
X^i
\]
be the continuous-time process which represents the average fitness of
the individuals in $X$. The rigorous results established in~\cite{YEC}
are as follows:
\begin{itemize}
\item The centered process $X^C$, in which individual $i$ has fitness
$X^{C,i} = X^i - \overline{X}$, is ergodic and has a stationary
distribution $\pi$.
\item If
\[
c_2 = \frac{1}{N}\sum_{i=1}^N
\bigl(X^{C,i}\bigr)^2
\]
is the variance of the centered process under the stationary
distribution, then
\[
E^\pi[\overline{X}_t] = \bigl(\mu(2q-1)+\gamma
E^\pi[c_2]\bigr)t,
\]
where $E^\pi$ means that the initial configuration of $X$ is chosen
according to the stationary distribution $\pi$.
\item For any $\delta> 0$ there exists $N_0$ large enough so that for
all $N \geq N_0$ we have $E\pi[\overline{X}_1] \geq\log^{1-\delta} N$.
\end{itemize}
It is difficult to say anything rigorous about $E^\pi[c_2]$ so other
methods are needed to compute $E[\overline{X}_t]$. The third result of
\cite{YEC} shows that if there is a positive ratio of beneficial
mutations then a large enough population will increase in fitness over
time. A paper by Etheridge and Yu~\cite{EY} provides further results
pertaining to this model.

Other similar models can be found in the biological literature. In
these models the density of the particles is assumed to act as a
traveling wave in time. The bulk of the wave behaves approximately
deterministically and the random noise comes from the most fit classes
of individuals. One tries to determine how quickly the fittest classes
advance and pull the wave forward. This traveling wave approach is used
in~\cite{YE} and~\cite{YEC} to approximate the rate of evolution as
$O(\log N/(\log\log N)^2)$. For other work in this direction see
Rouzine, Brunet and Wilke~\cite{RBW}, Brunet, Rouzine and Wilke \cite
{BRW}, Desai and Fisher~\cite{DF} and Park, Simon and Krug \cite
{PSK}. Using nonrigorous arguments, these authors get estimates of
$O(\log N)$, $O(\log N/\log\log N)$ and $O(\log N/(\log\log N)^2)$,
where the differences depend\vadjust{\goodbreak} on the details of the models that they
analyze. For more motivation and details concerning this model, please
see the Introduction in~\cite{YEC}.

Motivated by applications to cancer development, Durrett and Mayberry
have established rigorous results for a similar model in~\cite{DM}.
They consider two models in which all mutations are beneficial and the
mutation rate tends to 0 as the population size tends to infinity. In
one of their models the population size is fixed and in the other it is
exponentially increasing. For the model with the fixed population size
they show that the rate at which the average fitness is expected to
increase is $O(\log N)$. By considering the expected number of
individuals that have fitness $k$ at time $t$, they establish
rigorously that the density of the particles in their model will act as
a traveling wave in time.

Our result is the following theorem.
%
%
\begin{Theo} \label{Theorem}
Let $X_0^i = 0$ for $1 \leq i \leq N$. There exists a positive constant
$C$ which may depend on $\mu$, $q$ and $\gamma$ such that for $N$
large enough
\[
\frac{E[\overline{X}]}{t} \leq\frac{C\log N}{(\log\log N)^2}
\]
for all $t \geq\log\log N$.
\end{Theo}

A difference between the result in~\cite{YEC} and our result is that
in~\cite{YEC} the initial state of the process is randomly chosen
according to the stationary distribution $\pi$, while we make the
assumption that all of the individuals initially have fitness 0.

The statements of the propositions needed to prove Theorem \ref
{Theorem} and the proof of Theorem~\ref{Theorem} are included in
Section~\ref{sec2}. At the end of the paper there is a table which
includes the
notation that is used throughout the paper and the \hyperref
[app]{Appendix} that includes
some general results on branching processes.

\section{\texorpdfstring{Proof of Theorem \protect\ref{Theorem}}{Proof of Theorem 1}}\label{sec2}
Before stating the propositions we use to prove the theorem we need to
establish some notation. Let $X_t^+ = \max\{X_t^i\dvtx1 \leq i \leq N\}$
be the maximum fitness of any individual at time $t$ and $X_t^- = \min
\{X_t^i\dvtx1 \leq i \leq N\}$ be the minimum fitness of any
individual at
time $t$. Define the width of the process to be $W_t = X_t^+-X_t^-$ and
define $D_t = X_t^+ - X_0^+$ be the distance the front of the process
has traveled by time $t$. Theorem~\ref{Theorem} states that all
individuals initially have fitness~0. Therefore, a bound on $D_t$
immediately yields a bound on $\overline{X}_t$. The bounds we
establish on $D_t$ will depend on the width, $W_t$.

Let $w = w(N)$ be any positive, increasing function that satisfies
\[
\lim_{N\rightarrow\infty}w(N) = \infty
\quad\mbox{and}\quad\lim_{N\rightarrow
\infty}\frac{w(N)}{\log\log N}
= 0.
\]
Let $\mathcal{W} = \lfloor w\log N/\log\log N \rfloor$ and $\mathcal
{T} = w^{-1/2}\log\log N$. Heuristically, we conjecture that $W_t$ is
typically of size $O(\log N/\log\log N)$ so $\mathcal{W}$ is larger
than the typical width of $X$. With probability\vadjust{\goodbreak} tending to 1, selection
should cause any width larger than $\mathcal{W}$ to shrink within
$\mathcal{T}$ time units. Because the width is a stochastic process,
we are motivated to make the following definitions:
\begin{eqnarray*}
t_1 &=& 0,
\\[-2pt]
s_n &=& \inf\{t \geq t_n\dvtx W_t \geq2
\mathcal{W}\} \qquad\mbox{for } n \geq1,
\\[-2pt]
t_n &=& \inf\{t \geq s_{n-1}\dvtx W_t <
\mathcal{W}\} \qquad\mbox{for } n \geq2,
\\[-2pt]
Y_i &=& \sup_{s_i \leq t \leq t_{i+1}}D_t - D_{s_i}
\qquad\mbox{for }i \geq1,
\\[-2pt]
N_t &=& \max\{i\dvtx s_i \leq t\} \qquad\mbox{for }t \geq0.
\end{eqnarray*}
Note that $s_n$ and $t_n$ exist for all $n \geq1$ with probability 1.

We define branching processes $Z^{k,\uparrow}$ for $k \geq0$ which
have the following dynamics:
\begin{itemize}
\item Initially there are $N$ particles of type $k$ in $Z_0^{k,\uparrow}$.
\item Each particle changes from type $i$ to $i+1$ at rate $\mu$.
\item A particle of type $i$ branches at rate $\gamma i+1$ and, upon
branching, the new particle is also type $i$.
\end{itemize}
Let $\overline{M}{}^{k,\uparrow}_t$ be the maximum type of any particle
in $Z_t^{k,\uparrow}$ and let \mbox{$M_t^{k,\uparrow} = \overline
{M}{}^{k,\uparrow}_t-k$}, so that $M_0^{k,\uparrow} = 0$. Note that we
refer to individuals in branching processes as particles to distinguish
them from the individuals in $X$. This will make the coupling arguments
in the next section more clear.

We define a stochastic process $X'$ that will be coupled with $X$ as
described in the proof of Proposition~\ref{TheoProp1} for reasons that
will become clear shortly. Let $\{\mathcal{Z}^n\}_{n=0}^\infty$ be an
i.i.d. sequence of continuous-time stochastic processes which each have
the same distribution as $Z^{\mathcal{W},\uparrow}$. Let
$\overline{\mathcal{M}}{}^n_t$ be the maximum type of any particle in $\mathcal
{Z}_t^n$ and let $\mathcal{M}_t^n = \overline{\mathcal
{M}}{}^n_t-\mathcal{W}$ so that $\mathcal{M}_0^n = 0$ for all $n$. Define
\[
X_t' = \cases{ %
X_0^+ +
\mathcal{M}_t^0, &\quad if $t \in[0,\mathcal{T}]$,
\vspace*{1pt}\cr
X_{i\mathcal{T}}' + \mathcal{M}_{t-i\mathcal{T}}^i, &\quad
if $t \in\bigl(i\mathcal{T}, (i+1)\mathcal{T}\bigr]$\qquad for any integer $i \geq1$,}
\]
and $D_t' = X_t' - X_0^+$. The idea is that $D_t'$ is the maximum type
of any particle in a branching process $X'$ that has the same
distribution as $Z^{\mathcal{W},\uparrow}$ except that at each time
$i\mathcal{T}$ we restart the branching process so that there are once
again $N$ particles of type $\mathcal{W}$. For each integer $i \geq0$
at time $i\mathcal{T}$, the $N$ particles initially have type $D_t'$
which is the maximum type achieved by any particle in $X_t'$ up to time $t$.

Now we are able to state the four propositions used to prove Theorem
\ref{Theorem}. Proposition~\ref{TheoProp1} is a result of the
coupling of $X$ and $X'$.
%
%
\begin{Prop} \label{TheoProp1}
Let $X_0^i = 0$ for $1 \leq i \leq N$. Then
\[
D_t \leq D_t'+\sum
_{i=1}^{N_t}Y_i
\]
for all times $t \geq0$.\vadjust{\goodbreak}
\end{Prop}
%
%
\begin{Prop} \label{TheoProp4}
Let $X_0^i = 0$ for $1 \leq i \leq N$. For $N$ large enough we have
\[
\sup_{t \in[\mathcal{T},\infty)} \frac{E[D_t']}{t} \leq\frac
{6\mathcal{W}}{\mathcal{T}}.
\]
\end{Prop}

With the initial condition $X_0^i = 0$ for $1 \leq i \leq N$, we let
$\mathcal{F} = \{\mathcal{F}_t\}_{t \geq0}$ be the natural
filtration associated with $X$.
%
%
\begin{Prop} \label{TheoProp2}
Let $X_0^i = 0$ for $1 \leq i \leq N$. For $N$ large enough we have
$E[Y_i|\mathcal{F}_{s_i}] \leq5\mathcal{W}$ for all $i \geq1$.
\end{Prop}
%
%
\begin{Prop} \label{TheoProp3}
Let $X_0^i = 0$ for $1 \leq i \leq N$. For $N$ large enough,
\[
\sup_{s \in[0,\infty)}\frac{1}{s}E[N_s] \leq
\frac{1}{\mathcal{T}}.
\]
\end{Prop}
\begin{pf*}{Proof of Theorem~\ref{Theorem}}
Fix $t \geq\log\log N$. It follows by definition of $\mathcal{T}$
that $t > \mathcal{T}$ so that the hypotheses of the preceding four
propositions are satisfied. There exists $N_0$ which does not depend on
$t$ such that for any $N \geq N_0$ we have
\begin{eqnarray*}
E \biggl[\frac{D_t}{t} \biggr] & \leq & E \biggl[\frac{D_t'+\sum
_{i=1}^{N_t}Y_i}{t} \biggr]
\qquad\mbox{by Proposition~\ref{TheoProp1}}
\\
& = &E \biggl[\frac{D_t'}{t} \biggr] + E \biggl[\frac{\sum
_{i=1}^{N_t}Y_i}{t} \biggr]
\\
& \leq &\frac{6\mathcal{W}}{\mathcal{T}} + \frac{1}{t}E \Biggl[\sum
_{i=1}^{N_t}Y_i \Biggr] \qquad\mbox{by
Proposition~\ref{TheoProp4}}
\\
& = &\frac{6\mathcal{W}}{\mathcal{T}} + \frac{1}{t}\sum_{i=1}^\infty
E[Y_i 1_{\{N_t \geq i\}}]
\\
& = &\frac{6\mathcal{W}}{\mathcal{T}} + \frac{1}{t}\sum_{i=1}^\infty
E\bigl[E[Y_i 1_{\{N_t \geq i\}}|\mathcal{F}_{s_i}]\bigr]
\\
& = &\frac{6\mathcal{W}}{\mathcal{T}} + \frac{1}{t}\sum_{i=1}^\infty
E\bigl[1_{\{N_t \geq i\}}E[Y_i |\mathcal{F}_{s_i}]\bigr]
\\
& \leq &\frac{6\mathcal{W}}{\mathcal{T}} + \frac{5\mathcal
{W}}{t}\sum_{i=1}^\infty
E[1_{\{N_t \geq i\}}] \qquad\mbox{by Proposition~\ref{TheoProp2}}
\\
& = &\frac{6\mathcal{W}}{\mathcal{T}} + \frac{5\mathcal
{W}}{t}E[N_t]
\\
& \leq &\frac{6\mathcal{W}}{\mathcal{T}} + \frac{5\mathcal
{W}}{\mathcal{T}} \qquad\mbox{by Proposition
\ref{TheoProp3}}
\\
& = &\frac{11w^{1/2}\log N}{(\log\log N)^2}.
\end{eqnarray*}
Since $w$ may go to infinity arbitrarily slowly with $N$ there must
exist a constant $C$ such that
\[
\frac{E[D_t]}{t} \leq\frac{C\log N}{(\log\log N)^2}
\]
for all $t \geq\log\log N$. This immediately gives a bound on
$E[\overline{X}_t]/t$.
\end{pf*}

\section{Bounding the rate when the width is small}\label{sec3}
Through the use of branching processes we establish a bound on $D_t$
that depends on the width. We will make use of the strong Markov
property of $X$ at the times $s_n$ and $t_n$ for $n \geq1$. For this
reason, many of the statements we prove below will include conditions
for which $W_0 > 0$ even though according to the conditions of Theorem
\ref{Theorem} we have $W_0 = 0$. In this section we establish a small
upper bound for $D_t$ on the time intervals $[t_n,s_n)$.

The following proofs will involve coupling $X$ with various branching
processes. While the individuals in $X$ each have an integer value that
we refer to as the fitness of the individual, the particles in a
branching process will each be given an integer value that we refer to
as the type of the particle. Let $Z^C = \{Z_t^C\}_{t \geq0}$ be a
multi-type Yule process in which there are initially $N$ particles of
type 0. Particles increase from type $i$ to type $i+1$ at rate $\mu$
and branch at rate $C$. When a particle of type $i$ branches, the new
particle is also type $i$. Let $M_t^C$ be the maximum type of any
particle at time $t$.

The next proposition will give a lower bound on the fitness of any
individual up to time $t$ given that we know the least fitness at time
0 is $X_0^-$. We do this by establishing an upper bound on the amount
that any individual will decrease in fitness. Let
\[
S_t = \sup_{0 \leq s \leq t}\bigl(X_0^- -
X_s^-\bigr).
\]

%
%
\begin{figure}

\includegraphics{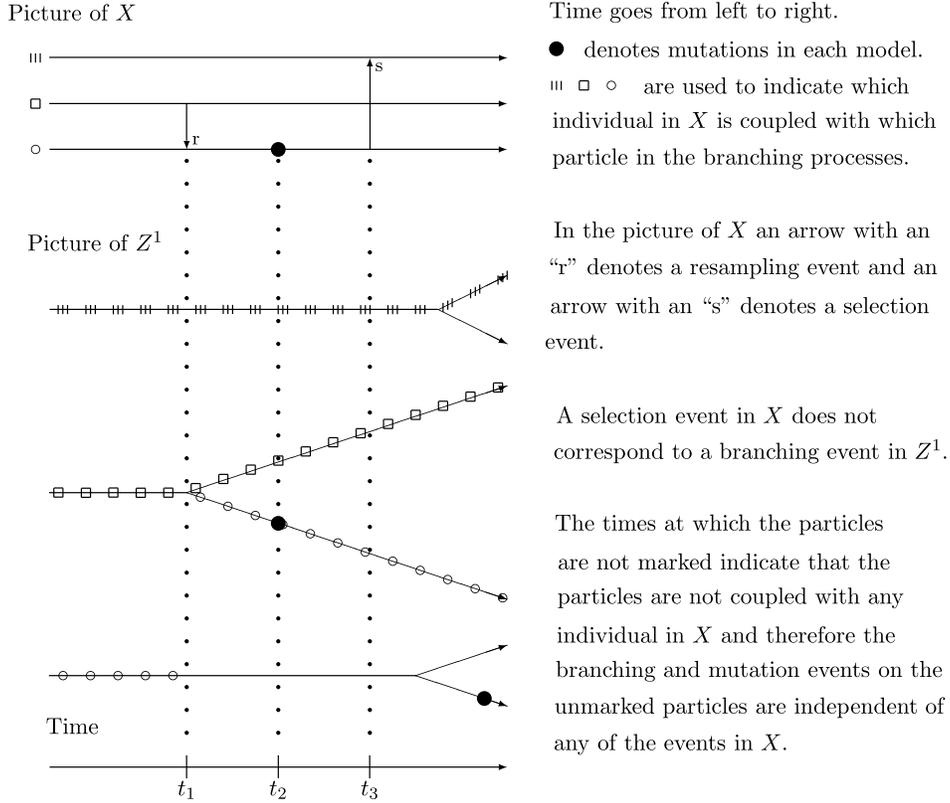}

\caption{Picture of the coupling of $X$ with $Z^1$ when $N = 3$.}\label{fig1}
\end{figure}

%
\begin{Prop} \label{BackSpeed}
For any population size $N$, initial configuration $X_0$, time $t \geq
0$ and natural number $l$,
\[
P(S_t \geq l) \leq\frac{N (t\mu)^l e^t}{l!}.
\]
\end{Prop}
\begin{pf}
By Lemma~\ref{YuleLem} in the \hyperref[app]{Appendix} we have
\[
P\bigl(M_t^1 \geq l\bigr) \leq\frac{N(t\mu)^le^t}{l!}
\]
for any population size $N$, time $t \geq0$ and natural number $l$.
Note that from our notation above $Z^1$\vadjust{\goodbreak} is a Yule process with
branching rate 1. To complete the proof we establish a
coupling\vspace*{1pt} between $X$ and $Z^1$ such that for any
population size $N$ and time $t \geq0$ we have $M_t^1 \geq S_t$. See
Figure~\ref{fig1} for an illustration of the coupling. At all times
every individual in $X$ will be paired with one particle in $Z^1$. The
coupling is as follows:
\begin{itemize}
\item We initially have a one-to-one pairing of each individual $i$ in
$X_0$ with each particle $i$ in $Z_0^1$.
\item The particle in $Z^1$ that is paired with individual $i$ will
increase in type by 1 only when individual $i$ gets a mutation.
\item For each individual $i$ in $X$ and each $j \neq i$, individual
$j$ is replaced by individual $i$ at rate $1/N$ due to resampling
events. If individual $i$ replaces individual $j$ due to resampling,
then the particle labeled $i$ in $Z^1$ branches. If particle $i$ has a
higher type than particle $j$, then the new particle is paired with
individual $j$. The particle that was paired with individual $j$ before
the branching event is no longer paired with any individual in $X$. If
particle $i$ has a lower type than particle $j$ then the particle that
was paired with individual $j$ remains paired with individual $j$ and
the new particle is not paired with any individual in $X$.
\item The particle paired with individual $i$ in $Z^1$ branches at rate
$1/N$ and these branching events are independent of any of the events
in $X$. When the particle paired with individual $i$ branches due to
these events, the new particle is not paired with any individual in $X$.
\item Any particles in $Z^1$ that are not paired with an individual in
$X$ branch and acquire mutations independently of $X$. The selection
events in $X$ are independent of any events in $Z^1$.
\end{itemize}

Let $R^i$ be the type of the particle in $Z^1$ that is paired with
individual $i$ and let
\[
S_s^i = \sup_{0\leq r \leq s}\bigl(X_0^--X_r^i
\bigr).
\]
To show $M_t^1 \geq S_t$ it is enough to show $R_t^i \geq S_t^i$ for
all $i$. Initially $S_0^i \leq R_0^i = 0$ for all $i$. Note that both
$s \mapsto S_s^i$ and $s \mapsto R_s^i$ are increasing functions and
that increases in these functions correspond to decreases in $X^i$.

When individual $i$ gets a mutation, $R^i$ increases by 1. However, if
individual $i$ gets a mutation at time $s$, then $S^i$ will only
increase by 1 if $S_{s-}^i = X_0^- - X_{s-}^i$ and the mutation is
deleterious. Therefore, if individual $i$ gets a mutation at time $s$
and $S_{s-}^i \leq R_{s-}^i$, then
\[
S_s^i \leq S_{s-}^i+1 \leq
R_{s^-}^i+1 = R_s^i.
\]

Suppose individual $j$ is replaced by individual $i$ due to a
resampling event at time $s$ and that both $S_{s-}^j \leq R_{s-}^j$ and
$S_{s-}^i \leq R_{s-}^i$ hold. With probability 1 we have $S_s^i =
S_{s-}^i$ and $R_s^i = R_{s-}^i$. If $X_0^- - X_s^i \leq S_{s-}^j$ then
$S_{s-}^j = S_s^j$. From this it follows that $S_s^j \leq R_s^j$. If
$X_0^- - X_s^i > S_{s-}^j$ then $S_s^j = X_0^--X_s^i \leq S_s^i \leq
R_s^i$. If $R_s^i \geq R_{s-}^j$, then by the definition of the
coupling, $R_s^j = R_s^i$. If $R_s^i < R_{s-}^j$, then by definition of
the coupling, $R_s^j = R_{s-}^j$. Therefore, $R_s^j \geq R_s^i$ which
gives us $S_s^j \leq R_s^j$.

Selection events will never increase $S^i$ and since $S^i$ and $R^i$
are increasing in time, a selection event at time $s$ will preserve the
inequality $S_s^i \leq R_s^i$. This shows that any event that occurs at
time $s$ which may change the fitness of an individual $i$ in $X$ will
preserve the inequality $S_s^i \leq R_s^i$. Since the result holds for
each individual~$i$, we have $S_t \leq M_t^1$.
\end{pf}

We now wish to bound the distance the front of the wave moves as a
function of the initial width.
%
%
\begin{Prop} \label{UpBound}
For any initial configuration $X_0$, fixed time $t \geq0$ and any
integer $l \geq0$, we have
\[
P\Bigl(\sup_{0\leq s \leq t} D_s > l\Bigr) \leq\frac{2N(t\mu
)^le^{(\gamma
(W_0+2l)+\mu+1)t}}{(l-1)!}.
\]
\end{Prop}
\begin{pf}
Recall that $W_0$ is the width of $X$ at time 0. We first establish a
coupling between $X$ and $Z^{W_0+k,\uparrow}$ for each integer $k \geq
0$. See Figure~\ref{fig2} for an illustration of the coupling.
%
%
\begin{figure}

\includegraphics{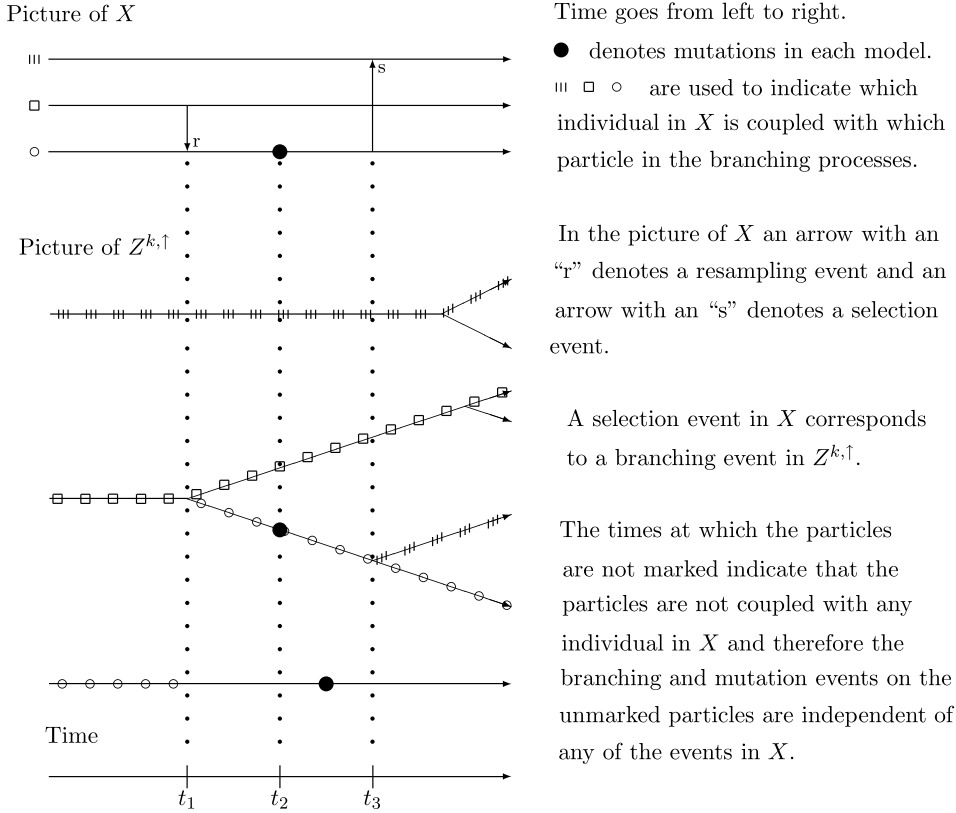}

\caption{Picture of the coupling of $X$ with $Z^{k,\uparrow}$ when $N
= 3$.}\label{fig2}
\end{figure}
Let $T^k= \inf\{
t\dvtx S_t > k\}$ for $k \geq1$. Every individual in $X$ will be paired
with one particle in $Z^{W_0+k,\uparrow}$ until time $T^k$. We couple
$Z^{W_0+k,\uparrow}$ with $X$ for all times $t \in[0,T^k)$ as
follows:
\begin{itemize}
\item We initially have a one-to-one pairing of each individual $i$ in
$X_0$ with each particle $i$ in $Z_0^{W_0+k,\uparrow}$. When a
particle in $Z_t^{W_0+k,\uparrow}$ is coupled with individual $i$, we
refer to the particle as particle $i$.
\item Particle $i$ increases in type by 1 only when individual $i$ gets
a mutation.
\item For each individual $i$ in $X$ and each $j \neq i$, individual
$j$ is replaced by individual $i$ at rate $1/N$ due to resampling
events. If individual $i$ replaces individual $j$ due to resampling,
then particle $i$ branches. If particle $i$ has a higher type than
particle $j$, then the new particle is paired with individual $j$. The
particle that was paired with individual $j$ before the branching event
is no longer paired with any individual in $X$. If particle $i$ has a
lower type than particle $j$, then the particle that was paired with
individual $j$ remains paired with individual $j$ and the new particle
is not paired with any individual in $X$.
\item Additionally, particle $i$ branches at rate $1/N$ and these
branching events are independent of any of the events in $X$. When
particle $i$ branches due to these events the new particle is not
paired with any individual in $X$.
\item In $X$ there is a time dependent rate $\gamma U_s^i$ at which
individuals $j \neq i$ are replaced by individual $i$ due to selection
events, namely,
\[
U_s^i = \frac{1}{N}\sum
_{j=1}^N \bigl(X_s^i -
X_s^j\bigr)^+.
\]
If individual $j$ is replaced by individual $i$ in $X$ due to a
selection event, then particle $i$ branches. If particle $i$ has a
higher type than particle $j$, then the new particle is paired with
individual $j$. The particle that was paired with individual $j$ before
the branching event is no longer paired with any individual in $X$. If
particle $i$ has a lower type than particle $j$, then the particle that
was paired with individual $j$ remains paired with individual $j$. The
new particle is not paired with any individual in $X$.
\item Additionally, particle $i$ branches at a time dependent rate
$\gamma(R_t^{i,k}-U_t^i)$ where $R_t^{i,k}$ is the type of particle
$i$. These branching events are independent of any of the events in
$X$. When such a branching event occurs, the new particle is not paired
with any individual in $X$.
\item Any particles in $Z^{W_0+k,\uparrow}$ that are not paired with
an individual in $X$ branch and change type independently of $X$.
\end{itemize}
Fix $k \geq1$. For the above coupling between $X$ and
$Z^{W_0+k,\uparrow}$ to be well defined until time $T^k$, we need
$R_t^{i,k}-U_t^i \geq0$ for all $i \in\{1,\ldots,N\}$ and for all
times $t \in[0,T^k)$. Let $\overline{T}{}^{k,i} = \inf\{t\dvtx R_t^{i,k} -
U_t^i < 0\}$. The coupling between $X$ and $Z^{W_0+k,\uparrow}$ is
well defined until time $\overline{T}{}^k = \min\{\overline{T}{}^{k,i}\dvtx1
\leq i \leq N\}$. We will show that $T^k \leq\overline{T}{}^k$.

Let
\[
\overline{S}{}^i_t = \sup_{0\leq s \leq t}
\bigl(X_s^i-X_0^+\bigr) \quad\mbox{and}\quad
\overline{R}{}^{i,k}_t = R_0^{i,k}-W_0-k.
\]
Initially $\overline{S}{}^i_0 \leq\overline{R}{}^{i,k}_0 = 0$ for all
$i$. Note that both $t \mapsto\overline{S}{}^i_t$ and $t \mapsto
R_t^{i,k}$ are increasing functions, from which it follows that $t
\mapsto\overline{R}{}^{i,k}_t$ is also an increasing function.

When individual $i$ gets a mutation, $\overline{R}{}^{i,k}$ increases by
1. However, if individual $i$ gets a mutation at time $s$ then
$\overline{S}{}^i$ will only increase by 1 if $\overline{S}{}^i_{s-} =
X_{s-}^i - X_0^+$ and the mutation is beneficial. Therefore, if
individual $i$ gets a mutation at time $s$ and $\overline{S}{}^i_{s-}
\leq\overline{R}{}^{i,k}_{s-}$, then
\[
\overline{S}{}^i_s \leq\overline{S}{}^i_{s-}+1
\leq\overline{R}{}^{i,k}_{s^-}+1 = \overline{R}{}^{i,k}_s.
\]

Suppose individual $j$ is replaced by individual $i$ due to a
resampling or selection event at time $s$ and that both $\overline
{S}{}^j_{s-} \leq\overline{R}{}^{j,k}_{s-}$ and $\overline{S}{}^i_s =
\overline{S}{}^i_{s-} \leq\overline{R}{}^{i,k}_{s-} = \overline
{R}{}^{i,k}_s$ hold. If $X_s^i - X_0^+ \leq\overline{S}{}^j_{s-}$, then
$\overline{S}{}^j_{s-} = \overline{S}{}^j_s$. It follows that $\overline
{S}{}^j_s \leq\overline{R}{}^{j,k}_s$. If $X_s^i-X_0^+ > \overline
{S}{}^j_{s-}$ then $\overline{S}{}^j_s = X_0^--X_s^i \leq\overline
{S}{}^i_s \leq\overline{R}{}^{i,k}_s$. If $\overline{R}{}^{i,k}_s \geq
\overline{R}{}^{j,k}_{s-}$, then by the definition of the coupling,
$\overline{R}{}^{j,k}_s = \overline{R}{}^{i,k}_s$. If $\overline
{R}{}^{i,k}_s < \overline{R}{}^{j,k}_{s-}$, then by definition of the
coupling, $\overline{R}{}^{j,k}_s = \overline{R}{}^{j,k}_{s-}$.
Therefore, $\overline{R}{}^{j,k}_s \geq\overline{R}{}^{i,k}_s$ which
gives us $\overline{S}{}^j_s \leq\overline{R}{}^{j,k}_s$.

For any time $s < T^k$ we have $R^{i,k}_s \geq\overline{S}{}^i_s+W_0+k
\geq X^i_s-X_0^++W_0+k = X^i_s-X_0^-+k$. If there were $N$ individuals
with fitness $X_0^--k$ at time $s \in[0,\overline{T}{}^{k,i})$, then
the rate at which individual $i$ replaces these $N$ individuals due to
selection is $\gamma(X^i_s-X_0^-+k)$. However, for any time $s < T^k$,
there are fewer than $N$ individuals being replaced by individual $i$
due to selection and they will all have fitnesses at least as large as
$X_0^--k$. This gives us a bound on the rate at which resampling events
occur on individual $i$ before time $T^k$, namely, $U^i_s \leq
X^i_s-X_0^-+k \leq R^{i,k}_s$ for all $s \in[0,T^k)$. This shows that
$T^k \leq\overline{T}{}^{k,i}$ for all $i$. Hence, $T^k \leq\overline
{T}{}^k$ and the coupling is well defined until time $T^k$.

We have shown that any event that occurs at time $s \in[0,T^k)$ which
may change the fitness of an individual $i$ in $X$ will preserve the
inequality $\overline{S}{}^i_s \leq\overline{R}{}^{i,k}_s$. Since the
result holds for each individual $i$, for any $s \in[0,T^k)$ we have
\[
\sup_{0 \leq r \leq s} D_r = \sup_{1 \leq i \leq N}
\overline{S}{}^i_s \leq\sup_{1 \leq i \leq N}
\overline{R}{}^{i,k}_s \leq M_s^{W_0+k,\uparrow}.
\]

Note that if $\sup_{0\leq s \leq t}(X_0^- - X_s^-) \leq k$ then $t <
T^k$. On the event\break $\{\sup_{0\leq s \leq t}(X_0^- - X_s^-) \leq k\}$
we have $M_t^{W_0+k,\uparrow} \geq\sup_{0\leq s \leq t}D_s$. This
allows us to do the following computation:
%
%
\begin{eqnarray}
\label{simpEq} P\Bigl(\sup_{0\leq s \leq t} D_s > l\Bigr) & = & \sum
_{i=0}^\infty P\Bigl(\Bigl\{
\sup_{0\leq s \leq t} D_s > l\Bigr\} \cap\Bigl\{\sup_{0\leq s \leq t}
\bigl(X_0^- - X_s^-\bigr) = i\Bigr\}\Bigr)
\nonumber
\\[-2pt]
& \leq &\sum_{i=0}^\infty P\Bigl(\bigl
\{M_t^{W_0+i,\uparrow} > l\bigr\} \cap\Bigl\{\sup_{0\leq s \leq t}
\bigl(X_0^- - X_s^-\bigr) = i\Bigr\}\Bigr)
\nonumber
\\[-2pt]
& \leq &\sum_{i=0}^\infty P\Bigl(\bigl
\{M_t^{W_0+i,\uparrow} > l\bigr\} \cap\Bigl\{\sup_{0\leq s \leq t}
\bigl(X_0^- - X_s^-\bigr) \geq i\Bigr\}\Bigr)
\nonumber
\\[-2pt]
& \leq &\sum_{i=0}^\infty P
\bigl(M_t^{W_0+i,\uparrow} > l\bigr) \wedge P\Bigl(\sup_{0\leq s \leq t}
\bigl(X_0^- - X_s^-\bigr) \geq i\Bigr)
\nonumber
\\[-2pt]
& \leq &\sum_{i=0}^\infty P
\bigl(M_t^{W_0+i,\uparrow} > l\bigr) \wedge\biggl(\frac{N (t\mu)^i e^t}{i!}
\biggr) \qquad\mbox{by Proposition~\ref{BackSpeed}}
\nonumber
\\[-2pt]
& \leq &\sum_{i=0}^\infty\biggl(
\frac{N (t\mu)^l e^{(\gamma
(W_0+i+l)+1)t}}{l!} \biggr) \wedge\biggl(\frac{N (t\mu)^i
e^t}{i!} \biggr) \nonumber\\[-2pt]
&&\mbox{by
Lemma~\ref{ZUpBound} in the \hyperref[app]{Appendix}}
\nonumber
\\[-2pt]
& \leq &\frac{N (t\mu)^l e^{(\gamma(W_0+l)+1)t}}{l!} \sum_{i=0}^{l-1}
e^{i\gamma t} + Ne^t \sum_{i=l}^\infty
\frac{(t\mu
)^i}{i!}
\nonumber
\\[-2pt]
& \leq &\frac{N (t\mu)^l e^{(\gamma(W_0+l)+1)t}}{l!}\cdot le^{l\gamma
t} + Ne^t \sum
_{i=l}^\infty\frac{(t\mu)^i}{i!}
\\[-2pt]
& \leq &\frac{N (t\mu)^l e^{(\gamma(W_0+2l)+1)t}}{(l-1)!} + \frac
{N(t\mu)^l e^{(\mu+1)t}}{l!} \nonumber\\[-2pt]
&&\mbox{by Lemma~\ref{TRTBoundLem}
in the \hyperref[app]{Appendix}}
\nonumber
\\[-2pt]
& \leq &\frac{2N(t\mu)^le^{(\gamma(W_0+2l)+\mu+1)t}}{(l-1)!}.\nonumber
\end{eqnarray}
\upqed\end{pf}

We now extend the bound we got on the least fit individuals in
Proposition~\ref{BackSpeed} to a slightly stronger result.
%
%
\begin{Dfn} \label{S}
Let $x \in\Z$ and let $\mathcal{S}_t^x \subset\{1,2,\ldots, N\}$
correspond to a collection of individuals at time $t$ which is
determined by the following dynamics:
\begin{itemize}
\item Initially, $\mathcal{S}_0^x$ consists of all individuals whose
fitness lies in the interval $(x,\infty)$.
\item If a resampling or selection event occurs at time $t$ and an
individual not in $\mathcal{S}_{t-}^x$ is replaced by a individual in
$\mathcal{S}_{t-}^x$, then it is added to $\mathcal{S}_t^x$.
\item If a beneficial mutation occurs at time $t$ on an individual not
in $\mathcal{S}_{t-}^x$ that causes its fitness to increase from $x$
to $x+1$, it is added to $\mathcal{S}_t^x$.
\item If a resampling event occurs at time $t$ to an individual in
$\mathcal{S}_{t-}^x$ and it is replaced by a individual not in
$\mathcal{S}_{t-}^x$, then it is removed from $\mathcal{S}_t^x$.
\end{itemize}
\end{Dfn}
Mutation and selection events do not cause individuals to be lost
from~$\mathcal{S}^x$. We now prove the following corollary to Proposition
\ref{UpBound}.
%
%
\begin{Cor} \label{BackSpeedLabel}
Let $A_t^{x,l}$ be the event that an individual in $\mathcal{S}_s^x$
has fitness in $(-\infty, x-l]$ for some time $s \in[0,t]$. For any
initial configuration $X_0$, time $t \geq0$ and any integer $l$,
\[
P\bigl(A_t^{x,l}\bigr) \leq\frac{2N(t\mu)^le^{(\gamma(W_0+2l)+\mu+1)t}}{(l-1)!}.
\]
\end{Cor}

Note that we cannot use the bound found in Proposition~\ref{BackSpeed}
because individuals not in $\mathcal{S}_t^x$ may move\vadjust{\goodbreak} to $\mathcal
{S}_t^x$ due to selection events. In the proof of Proposition \ref
{BackSpeed} the number of individuals with the least fitness cannot
increase due to selection events. However, the number of individuals
with the least fitness in $\mathcal{S}_t^x$ may increase due to
selection events involving individuals not in $\mathcal{S}_t^x$.
\begin{pf*}{Proof of Corollary~\ref{BackSpeedLabel}}
For $k \geq1$ let $X$ be coupled with $Z^{W_0+k,\uparrow}$ as in the
proof of Proposition~\ref{UpBound}. Let $T^k$, $R^{i,k}_t$ and
$\overline{R}{}^{i,k}_t$ be defined as they were in the proof of
Proposition~\ref{UpBound}. Define $\overline{T}{}^i_s = \{r \in[0,s]\dvtx
i \in\mathcal{S}_r^x\}$ and let
\[
S^i_s = \cases{ %
\displaystyle \sup_{r \in\overline{T}{}^i_s}
\bigl(x-X_r^i\bigr), &\quad if $\overline{T}{}^i_s
\neq\varnothing$,
\vspace*{2pt}\cr
-\infty, &\quad if $\overline{T}{}^i_s =
\varnothing$.}
\]

The goal is to show that for all $s \in[0,T^k)$ we have
\[
\sup_{1 \leq i \leq N} S^i_s \leq\sup_{1 \leq i \leq N}
\overline{R}{}^{i,k}_s \leq M_s^{W_0+k,\uparrow}.
\]
Note that we can only consider the coupling of $X$ with
$Z^{W_0+k,\uparrow}$ until time $T^k$ because after this time the
coupling is not well defined.

Initially all of the individuals in $\mathcal{S}_0^x$ have fitness in
$(x,\infty)$. Therefore, if $i \in\mathcal{S}_0^x$ then $S_0^i \leq
0 = \overline{R}{}^{i,k}_0$. If $i \notin\mathcal{S}_0^x$ then $S_0^i
= -\infty< \overline{R}{}^{i,k}_0$.

Suppose individual $i$ gets a mutation at time $s$ and that for any
time $s' \in[0,s-)$ we have $S_{s'}^i \leq\overline{R}{}^{i,k}_{s'}$.
Then $\overline{R}{}^{i,k}$ increases by 1. If $i \in\overline
{S}{}^x_{s-}$ then $S^i_s$ will only increase by 1 if $S_{s-}^i =
x-X^i_s$ and the mutation is deleterious. If $i \notin\mathcal
{S}_{s-}^x$ and the mutation does not cause the fitness of individual
$i$ to change from $x$ to $x+1$, then $S^i_s = S_{s-}^i$. If $i \notin
\mathcal{S}_{s-}^x$ and the mutation does cause the fitness of
individual $i$ to change from $x$ to $x+1$, then $S^i_s = S_{s-}^i \vee
0$. In any of these three cases, $S^i_s \leq\overline{R}{}^{i,k}_s$.

Suppose individual $j$ is replaced by individual $i$ due to a
resampling or selection event at time $s$ and that $S_{s-}^j \leq
\overline{R}{}^{j,k}_{s-}$ and $S_{s-}^i \leq\overline{R}{}^{i,k}_{s-}$.
If $i \notin\overline{S}{}^x_{s-}$ then $S_{s-}^j = S_s^j \leq
\overline{R}{}^{j,k}_{s-}$. Suppose $i \in\mathcal{S}_{s-}^x$. If $x -
X^i_s \leq S_{s-}^j$ then $S_{s-}^j = S_s^j$. From this it follows that
$S_s^j \leq\overline{R}{}^j_s$. If $x - X^i_s > S_{s-}^j$, then $S_s^j
= x-X^i_s \leq S^i_s \leq\overline{R}{}^i_s$. If $\overline{R}{}^i_s
\geq\overline{R}{}^j_{s-}$, then by the definition of the coupling,
$\overline{R}{}^j_s = \overline{R}{}^i_s$. If $\overline{R}{}^i_s <
\overline{R}{}^j_{s-}$, then by definition of the coupling, $\overline
{R}{}^j_s = \overline{R}{}^j_{s-}$. Therefore, $\overline{R}{}^j_s \geq
\overline{R}{}^i_s$ which gives us $S^j_s \leq\overline{R}{}^j_s$.

Note that if $\sup_{0\leq s \leq t}(X_0^- - X_s^-) \leq k$ then $t <
T^k$. Therefore, on the event $\{\sup_{0\leq s \leq t}(X_0^- - X_s^-)
\leq k\}$ we have $M^{W_0+k,\uparrow}_t \geq\sup_{1 \leq i \leq N}
S^i_s$. This allows us to do the following computation:
\begin{eqnarray*}
&&
P\Bigl(\sup_{0\leq s \leq t} \sup_{1 \leq i \leq N} S^i_s
> l\Bigr) \\
&&\qquad = \sum_{i=0}^\infty P\Bigl(\Bigl
\{\sup_{0\leq s \leq t} \sup_{1 \leq i \leq N} S^i_s > l
\Bigr\} \cap\Bigl\{\sup_{0\leq s \leq t}\bigl(X_0^- -
X_s^-\bigr) = i\Bigr\}\Bigr)
\\
&&\qquad \leq \sum_{i=0}^\infty P\Bigl(\bigl
\{M_t^{W_0+i,\uparrow} > l\bigr\} \cap\Bigl\{\sup_{0\leq s \leq t}
\bigl(X_0^- - X_s^-\bigr) = i\Bigr\}\Bigr).
\end{eqnarray*}
This is the same bound as equation (\ref{simpEq}) in the proof of
Proposition~\ref{UpBound}. Therefore, we have established the same bound.
\end{pf*}
\begin{pf*}{Proof of Proposition~\ref{TheoProp4}}
By definition $D_\mathcal{T}'$ has the same distribution as
$M_\mathcal{T}^{\mathcal{W},\uparrow}$ so by Lemma~\ref{ZUpBound}
in the \hyperref[app]{Appendix} we have
\[
P\bigl(D_\mathcal{T}' > l\bigr) \leq\frac{N(\mathcal{T}\mu)^l e^{(\gamma
(\mathcal{W}+l)+1)\mathcal{T}}}{l!}.
\]
Then
%
%
\begin{eqnarray}
\label{eq1} \frac{E[D_\mathcal{T}']}{2\mathcal{W}} & = & \frac
{1}{2\mathcal
{W}}\sum
_{l=0}^\infty P\bigl(D_\mathcal{T}'
> l\bigr)
\nonumber\\[-8pt]\\[-8pt]
& \leq &\frac{1}{2\mathcal{W}} \Biggl[2\mathcal{W} + \sum
_{l=2\mathcal{W}}^\infty\frac{N(\mathcal{T}\mu)^l e^{(\gamma
(\mathcal{W}+l)+1)\mathcal{T}}}{l!} \Biggr].\nonumber
\end{eqnarray}
By Lemma~\ref{TRTBoundLem} in the \hyperref[app]{Appendix} we have
%
%
\begin{equation}
\label{eq2} \sum_{l=2\mathcal{W}}^\infty
\frac{N(\mathcal{T}\mu)^l e^{(\gamma
(\mathcal{W}+l)+1)\mathcal{T}}}{l!} \leq\frac{N e^{(\gamma\mathcal
{W} + 1)\mathcal{T}}(\mathcal{T}\mu e^{\gamma\mathcal
{T}})^{2\mathcal{W}}e^{\mathcal{T} \mu e^{\gamma\mathcal
{T}}}}{(2\mathcal{W})!}.
\end{equation}

Note that for any $k \geq2$ both $D_{k\mathcal{T}}'-D_{(k-1)\mathcal
{T}}'$ and $D_\mathcal{T}'$ have the same distribution, namely, that
of $M_\mathcal{T}^{\mathcal{W}}$. Choose $t \in[k\mathcal{T},
(k+1)\mathcal{T})$ for some $k \geq1$. Because $D_t'$ is increasing
in $t$ we have
\[
\frac{D_t'}{t} \leq\frac{1}{k\mathcal{T}} \bigl(D_{(k+1)\mathcal
{T}}'
- D_{k\mathcal{T}}'+D_{k\mathcal{T}}'-
\cdots+D_{2\mathcal
{T}}'-D_\mathcal{T}'+D_\mathcal{T}'
\bigr).
\]
Therefore,
\[
\frac{E[D_t']}{t} \leq\frac{(k+1)E[D_\mathcal{T}']}{k\mathcal{T}} \leq
\frac{2E[D_\mathcal{T}']}{\mathcal{T}}.
\]

Let $t > \mathcal{T}$. Dividing both sides by $2\mathcal{W}/\mathcal
{T}$ and using the bounds found in equations (\ref{eq1}) and (\ref
{eq2}) gives us
\[
\frac{\mathcal{T}E[D_t']}{2t\mathcal{W}} \leq\frac{2E[D_\mathcal
{T}']}{2\mathcal{W}} \leq2 + \frac{N e^{(\gamma\mathcal{W} +
1)\mathcal{T}}(\mathcal{T}\mu e^{\gamma\mathcal{T}})^{2\mathcal
{W}}e^{\mathcal{T} \mu e^{\gamma\mathcal{T}}}}{2\mathcal
{W}(2\mathcal{W})!}.
\]
By Stirling's formula we have
\[
\frac{N e^{(\gamma\mathcal{W} + 1)\mathcal{T}}(\mathcal{T}\mu
e^{\gamma\mathcal{T}})^{2\mathcal{W}}e^{\mathcal{T} \mu e^{\gamma
\mathcal{T}}}}{2\mathcal{W}(2\mathcal{W})!} \sim\frac{N e^{(\gamma
\mathcal{W} + 1)\mathcal{T}}(\mathcal{T}\mu e^{\gamma\mathcal
{T}})^{2\mathcal{W}}e^{\mathcal{T} \mu e^{\gamma\mathcal
{T}}+2\mathcal{W}}}{(2\mathcal{W})^{2\mathcal{W}+1}\sqrt{4\pi
\mathcal{W}}} = e^x,
\]
where
\begin{eqnarray*}
x &=& \log N + \mathcal{T}\bigl(\gamma\mathcal{W}+1+\mu e^{\gamma\mathcal
{T}}\bigr)+2
\mathcal{W}\bigl(\log\bigl(\mathcal{T}\mu e^{\gamma\mathcal
{T}}\bigr)+1\bigr)\\
&&{}-(2
\mathcal{W}+1)\log(2\mathcal{W})-\log(4\pi\mathcal{W})/2.
\end{eqnarray*}
As $N \rightarrow\infty$ we have $x \sim-(2\mathcal{W}+1)\log
(2\mathcal{W}) \sim-2w\log N$. Therefore,
\[
\frac{\mathcal{T}E[D_t']}{2t\mathcal{W}} \leq3
\]
for $N$ large enough.
\end{pf*}
\begin{pf*}{Proof of Proposition~\ref{TheoProp1}}
We now couple $X$ with $X'$ by coupling $X$ with the sequence of
processes $\{\mathcal{Z}^m\}_{m=0}^\infty$. Let
\[
I_m = \bigl(m\mathcal{T},(m+1)\mathcal{T}\bigr] \cap\bigcup
_{n=1}^\infty[t_n,s_n)
\quad
\mbox{and}
\quad
J_m = (0,\mathcal{T}] \cap\bigcup
_{n=1}^\infty[t_n-m\mathcal{T},s_n-m
\mathcal{T}).
\]
For any $m \geq0$ we couple $X$ and $\mathcal{Z}^m$ as follows:
\begin{itemize}
\item The particles in $\mathcal{Z}_0^m$ are labeled
$1,2,\ldots,N$.\vspace*{1pt}
\item For any time in $I_m^C$ the process $X$ behaves independently of
$\mathcal{Z}^m$. For any time in $J_m^C$ the process $\mathcal{Z}^m$
behaves independently of the process $X$. During the time~$J_m^C$, if a
particle labeled $i$ in $\mathcal{Z}^m$ branches, the particle remains
labeled $i$ and the new particle is unlabeled.
\item The particle in $\mathcal{Z}^m$ that is paired with individual
$i$ will increase in type by 1 at time $t \in J_m$ only when individual
$i$ gets a mutation at time \mbox{$t+m\mathcal{T} \in I_m$}.
\item For each individual $i$ in $X$ and each $j \neq i$, individual
$j$ is replaced by individual $i$ at rate $1/N$ due to resampling
events. If individual $i$ replaces individual $j$ due to resampling at
time $t \in I_m$, then the particle labeled $i$ in $\mathcal{Z}^m$
branches at time $t-m\mathcal{T} \in J_m$. If particle $i$ has a
higher type than particle $j$, then the new particle is paired with
individual $j$. The particle that was paired with individual $j$ before
the branching event is no longer paired with any individual in $X$. If
particle $i$ has a lower type than particle~$j$, then the particle that
was paired with individual $j$ remains paired with individual $j$ and
the new particle is not paired with any individual in $X$.
\item The particle paired with individual $i$ in $\mathcal{Z}^m$
branches at rate $1/N$ for all times $t \in J_m$ and these branching
events are independent of any of the events in $X$. When the particle
paired with individual $i$ branches due to these events the new
particle is not paired with any individual in $X$.
\item In $X$ there is a time dependent rate $\gamma U^i_s$ at which
individuals $j \neq i$ are replaced by individual $i$ due to selection
events. If individual $j$ is replaced by individual $i$ in $X$ due to a
selection event at time $t \in I_m$, then the particle labeled\vadjust{\goodbreak} $i$ in
$\mathcal{Z}^m$ splits at time $t-m\mathcal{T} \in J_m$. If particle
$i$ has a higher type than particle $j$, then the new particle is
paired with individual~$j$. The particle that was paired with
individual $j$ before the branching event is no longer paired with any
individual in $X$. If particle $i$ has a lower type than particle $j$,
then the particle that was paired with individual $j$ remains paired
with individual $j$. The new particle is not paired with any individual
in $X$.
\item A particle labeled $i$ in $\mathcal{Z}^m$ splits at a
time-dependent rate $\gamma(R^{i,k}_t-U_t^i)$ for all times $t \in
J_m$ where $R^{i,k}_t$ is the type of particle $i$. These branching
events are independent of any of the events in $X$. When such a
branching event occurs, the new particle is not paired with any
individual in $X$.
\item Any particles in $\mathcal{Z}^m$ that are not paired with an
individual in $X$ branch and acquire mutations independently of $X$.
\end{itemize}

Observe the following bound for $D_t$:
\begin{eqnarray*}
D_t &\leq&\sum_{i=1}^{N_t-1}(D_{t_{i+1}}-D_{s_i})
+ \sum_{i=1}^{N_t}(D_{s_i}-D_{t_i})
+ \sup_{s_{N_t} \leq s \leq t_{N_t+1}} (D_s - D_{s_{N_t}})\\
&&{} +
\sup_{t_{N_t+1} \leq s \leq t} (D_s - D_{t_{N_t+1}}),
\end{eqnarray*}
where we consider the supremum over the empty set to be 0. By
definition we have
\[
\sum_{i=1}^{N_t-1}(D_{t_{i+1}}-D_{s_i})
+ \sup_{s_{N_t} \leq s \leq
t_{N_t+1}} (D_s - D_{s_{N_t}}) \leq\sum
_{i=1}^{N_t} Y_i.
\]

To finish the proof we will show
\[
\sum_{i=1}^{N_t}\sup_{t_i \leq s \leq s_i}(D_s-D_{t_i})
+ \sup_{t_{N_t+1} \leq s \leq t} (D_s - D_{t_{N_t+1}}) \leq
D_t'.
\]
To do this we define
\[
M_t = \sum_{i=1}^{N_t}
\sup_{t_i \leq s \leq s_i}(D_s-D_{t_i}) + \sup_{t_{N_t+1} \leq s \leq t}
(D_s - D_{t_{N_t+1}})
\]
for all times $t \geq0$. Suppose $M_s \leq D_s'$ for all $s \in[0,t)$
and a mutation, resampling or selection event occurs in $X$ at time
$t$. If $t \in(s_i, t_{i+1})$ for some $i \geq0$, then $M_{t-} = M_t$
because the process $M$ does not change on these time intervals. It is
possible that $D_t'$ changes, but $D_t'$ can only increase. Therefore,
$D_t' \geq M_t$. If $t \in[t_i, s_i] \cap(m\mathcal{T},(m+1)\mathcal
{T}]$ for some $i \geq0$ and $m \geq0$, then at time $t$ the
processes $X$ and $X'$ are coupled. More precisely, $X$ and $\mathcal
{Z}^m$ are coupled and the coupling has the same dynamics as the
coupling in Proposition~\ref{UpBound} except the time shift. The same
argument used in Proposition~\ref{UpBound} shows that $D_t' \geq M_t$
whether the individual changed fitness due to mutation, resampling or
selection. Since this inequality is preserved on any event that may
change $M_t$, it is true for all times $t$.
\end{pf*}

\section{Bounding the rate when the width is large}\label{sec4}
We consider what happens when the width is large in this section. By
large width we mean $W_t \geq\mathcal{W}$. The statements in this
section are easier to make when we consider an initial configuration of
$X$ such that $W_0 \geq\mathcal{W}$. Although the conditions of
Theorem~\ref{Theorem} state that $W_0 = 0$, we can wait for a random
time $\tau$ so that $W_\tau\geq\mathcal{W}$ and apply the strong
Markov property.

We begin this section by showing that when the width is large enough
the selection mechanism will cause the width to decrease quickly. We
give a labeling to the individuals that will help
us in this regard. Define the following subsets of $\R$:
\begin{eqnarray*}
I_1 &=& \bigl(-\infty, X_0^+-\tfrac{3}{16}W_0\bigr],
\\
I_2 &=& \bigl(X_0^+-\tfrac{3}{16}W_0,
X_0^+-\tfrac{2}{16}W_0\bigr],
\\
I_3 &=& \bigl(X_0^+-\tfrac{2}{16}W_0,
X_0^+-\tfrac{1}{16}W_0\bigr],
\\
I_4 &=& \bigl(X_0^+-\tfrac{1}{16}W_0,
\infty\bigr).
\end{eqnarray*}

We will label each individual in $X_0$ with two labels. For the first
labeling, we use $\mathfrak{a}$ to label the individuals in $I_1 \cup
I_2$, we use $\mathfrak{b}$ to label the individuals in $I_3$ and we
use $\mathfrak{c}$ to label the individuals in $I_4$. For the second
labeling we use $\mathfrak{a}'$ to label the individuals in $I_1$, we
use $\mathfrak{b}'$ to label the individuals in $I_2$ and we use
$\mathfrak{c}'$ to label the individuals in $I_3 \cup I_4$.

Let $\mathfrak{A}_t$, $\mathfrak{B}_t$ and $\mathfrak{C}_t$ denote
the number of individuals labeled $\mathfrak{a}$, $\mathfrak{b}$ and
$\mathfrak{c}$ at time~$t$, respectively. Let $\mathfrak{A}_t'$,
$\mathfrak{B}_t'$ and $\mathfrak{C}_t'$ denote the number of
individuals labeled $\mathfrak{a}'$, $\mathfrak{b}'$ and $\mathfrak
{c}'$ at time $t$, respectively.

The individuals change labels over time according to the following dynamics:
\begin{itemize}
\item Mutations: If the fitness of an individual labeled $\mathfrak
{a}$ increases so that it is in~$I_3$, then the individual is relabeled
$\mathfrak{b}$. If the fitness of a individual labeled $\mathfrak
{a}'$ increases so that it is in $I_2$, then the individual is
relabeled $\mathfrak{b}'$. Likewise, if the fitness of a individual
labeled $\mathfrak{b}$ increases so that it is in $I_4$, then it is
relabeled $\mathfrak{c}$ and if the fitness of a individual labeled
$\mathfrak{b}'$ increases so that it is in $I_3$, then it is relabeled
$\mathfrak{c}'$. Deleterious mutations do not cause individuals to be
relabeled.
\item Resampling: Any resampling event in which individual $i$ is
replaced by individual $j$ causes individual $i$ to inherit the labels
of individual $j$.
\item Selection: If an individual labeled $\mathfrak{a}$ is replaced
due to a selection event, it inherits the corresponding label of the
individual that replaced it. If an individual labeled $\mathfrak{a}'$
is replaced due to a selection event, it inherits the corresponding
label of the individual that replaced it. If an individual labeled
$\mathfrak{b}$ is replaced by an individual labeled $\mathfrak{c}$
due to a selection event, then the individual that was labeled
$\mathfrak{b}$ is relabeled $\mathfrak{c}$. If an individual labeled
$\mathfrak{b}'$ is replaced by\vadjust{\goodbreak} an individual labeled $\mathfrak{c}'$
due to a selection event, then the individual that was labeled
$\mathfrak{b}'$ is relabeled~$\mathfrak{c}'$. Any other selection
events do not cause the labels of the individuals to be changed.
\end{itemize}

Let $A_1$ be the event that there is an individual labeled $\mathfrak
{b}$ with fitness in $(-\infty, X_0^+-\frac{5}{32}W_0)$ for some time
$t \in[0,\mathcal{T}]$. Let $A_2$ be the event that there is an
individual labeled $\mathfrak{c}$ with fitness in $(-\infty,X_0^+ -
\frac{3}{32}W_0)$ for some time $t \in[0,\mathcal{T}]$. Let $A_1'$
be the event that there is an individual labeled $\mathfrak{b}'$ with
fitness in $(-\infty, X_0^+-\frac{7}{32}W_0)$ for some time $t \in
[0,\mathcal{T}]$. Let $A_2'$ be the event that there is an individual
labeled $\mathfrak{c}'$ with fitness in $(-\infty,X_0^+ - \frac
{5}{32}W_0)$ for some time $t \in[0,\mathcal{T}]$.
%
%
\begin{Lemma} \label{AALem}
Suppose $W_0 \geq\mathcal{W}$ for all $N$. Then
\[
P\bigl(A_1 \cup A_2 \cup A_1'
\cup A_2'\bigr) \rightarrow0 \qquad\mbox{as } N \rightarrow
\infty.
\]
\end{Lemma}
\begin{pf}
First we show the result for $A_1$. We apply Corollary \ref
{BackSpeedLabel} with $x = X_0^+-2W_0/16$, $t = t_0$ and $l = W_0/32$.
Recall that we had defined $\mathcal{S}_t^x$ in Definition~\ref{S}.
Because $x = X_0^+-2W_0/16$, we have that $\mathcal{S}_0^x$ consists
of all the individuals labeled $\mathfrak{b}$ or $\mathfrak{c}$.
Setting $t = \mathcal{T}$ and $l = W_0/32$ will make $A_t^{x,l}$ the
event that an individual labeled $\mathfrak{b}$ or $\mathfrak{c}$ has
fitness less than $X_0^+-\frac{5}{32}W_0$ by time $\mathcal{T}$. Note
that according to the relabeling dynamics, individual $i$ being labeled
$\mathfrak{b}$ or $\mathfrak{c}$ is equivalent to $i \in\mathcal
{S}^x$. Therefore, $A_1 \subset A_t^{x,l}$ and we get
\[
P(A_1) \leq P\bigl(A_t^l\bigr) \leq
\frac{2N(t\mu)^le^{(\gamma(W_0+2l)+\mu
+1)t}}{\lfloor l-1 \rfloor!}.
\]
Applying Stirling's formula we have
\[
\frac{2N(t\mu)^le^{(\gamma(W_0+2l)+\mu+1)t}}{\lfloor l-1 \rfloor!} \sim
\frac{2N(t\mu)^le^{(\gamma(W_0+2l)+\mu+1)t+\lfloor l-1 \rfloor
}}{\lfloor l-1 \rfloor^{\lfloor l-1 \rfloor}\sqrt{2\pi\lfloor l-1
\rfloor}} = e^x,
\]
where
\begin{eqnarray*}
x &=& \log(2N)+l\log(t\mu)+\bigl(\gamma(W_0+2l)+\mu+1\bigr)t+\lfloor
l-1 \rfloor\\
&&{}-\lfloor l-1 \rfloor\log\bigl(\lfloor l-1 \rfloor\bigr)-\log
\bigl(2\pi
\lfloor l-1 \rfloor\bigr)/2.
\end{eqnarray*}
As $N \rightarrow\infty$ we have $x \sim-\lfloor l-1 \rfloor\log
(\lfloor l-1 \rfloor) \sim-w\log N/32$. Therefore,
\[
P(A_1) \rightarrow0 \qquad\mbox{as }N \rightarrow\infty.
\]

We can apply Corollary~\ref{BackSpeedLabel} with $x = X_0^+-W_0/16$,
$t = \mathcal{T}$ and $l = W_0/32$ to get the same bound for $P(A_2)$.
By choosing $x$, $t$ and $l$ in this way, the event $A_t^{x,l}$ is the
event that an individual labeled $\mathfrak{c}$ has fitness less than
$X^+(0)-\frac{3}{32}W_0$ by time~$\mathcal{T}$. This shows that
$P(A_2)$ also tends to 0 as $N$ tends to infinity.

Likewise, to show $P(A_1')$ tends to 0 as $N$ goes to infinity we can
apply Corollary~\ref{BackSpeedLabel} with $x = X_0^+-\frac
{3}{16}W_0$, $t = \mathcal{T}$ and $l = W_0/32$,\vadjust{\goodbreak} and to show $P(A_2')$
tends to 0 as $N$ goes to infinity we can apply Corollary \ref
{BackSpeedLabel} with $x = X_0^+-\frac{2}{16}W_0$, $t = \mathcal{T}$
and $l = W_0/32$.
\end{pf}
%
%
\begin{Lemma} \label{CLem}
Suppose $W_0 \geq\mathcal{W}$ for all $N$. Let $T$ be a stopping time
whose definition may depend on $N$ such that $\mathfrak{C}_T' \geq
N/4$ for all $N$. Let $B_T = \inf\{t \geq T\dvtx X_t^- > X_0^+ - W_0/4\}
$. Then
\[
P\bigl(B_T1_{\{T<\mathcal{T}/2\}} > \tfrac{1}{2}\mathcal
{T}\bigr)
\rightarrow0 \qquad\mbox{as } N \rightarrow\infty.
\]
\end{Lemma}
\begin{pf}
Let $A_3'$ be the event that $\mathfrak{C}_t' \geq N/5$ for all times
$t \in[T,T+\frac{1}{2}\mathcal{T})$. The only way for an individual
labeled $\mathfrak{c}'$ to change its label is for it to be replaced
by an individual labeled $\mathfrak{a}'$ or $\mathfrak{b}'$ via a
resampling event. The rate at which individuals marked $\mathfrak{c}'$
undergo resampling events with individuals marked $\mathfrak{a}'$ or
$\mathfrak{b}'$ at time $t$ is
\[
\frac{\mathfrak{C}_t'(N-\mathfrak{C}_t')}{N} \leq\frac{N}{4}.
\]

Let $\{U_n\}_{n=0}^\infty$ be a simple random walk with $U_0 = N/4
\leq\mathfrak{C}_T'$. Let $T \leq t_1 < t_2 < \cdots$ be the times at
which individuals labeled $\mathfrak{c}'$ are involved in resampling
events with individuals that are not labeled $\mathfrak{c}'$ after
time $T$. We couple $\{U_n\}_{n=0}^\infty$ with $X$ so that if at time
$t_n$ an individual is labeled $\mathfrak{c}'$ due to a resampling
event, then $U_n = U_{n-1}+1$. If at time $t_n$ an individual loses the
label $\mathfrak{c}'$ due to a resampling event, then $U_n =
U_{n-1}-1$. To have $U_m < N/5$ for some $m$ satisfying $0 \leq m \leq
n$ we will need ${\max_{0 \leq m \leq n}}|U_m-U_0| \geq N/20$. It
follows from the reflection principle that there exists a constant $C$
such that $E[{\max_{0 \leq m \leq n}}|U_m-U_0|] \leq C\sqrt{n}$ for all
$n \geq0$. By Markov's inequality,
\[
P \Bigl(\max_{0 \leq m \leq n}|U_m-U_0| \geq N/20 \Bigr)
\leq C\sqrt{n}/N
\]
for some constant $C$.

Let $R$ be the number of resampling events that occur in the time
interval $[T,T+\frac{1}{2}\mathcal{T})$ that involve pairs of
individuals such that one is labeled $\mathfrak{c}'$ and the other is
not. Using Lemma~\ref{TRTBoundLem} in the \hyperref[app]{Appendix} and
the fact that the
rate at which resampling events occur is bounded above by $N/4$, we have
\[
P(R > k) \leq\sum_{i=k+1}^\infty
\frac{(N\mathcal{T})^i
e^{-N\mathcal{T}/8}}{8^i i!} \leq\frac{(N\mathcal{T})^k}{8^k k!}.
\]
Then
\begin{eqnarray*}
P\bigl(\bigl(A_3'\bigr)^C\bigr) & \leq & P
\Bigl(\Bigl\{\max_{0 \leq m \leq R}|U_m-U_0| \geq N/20\Bigr
\} \cap\bigl\{R \leq N^{3/2}\bigr\}\Bigr)
\\
&&{} + P\Bigl(\Bigl\{\max_{0 \leq m \leq R}|U_m-U_0| \geq
N/20\Bigr\} \cap\bigl\{R > N^{3/2}\bigr\}\Bigr)
\\
& \leq & P\Bigl(\Bigl\{\max_{0 \leq m \leq N^{3/2}}|U_m-U_0| \geq
N/20\Bigr\}\Bigr) + P\bigl(R > N^{3/2}\bigr)
\\[-2pt]
& \leq &\frac{C}{N^{1/4}} + \frac{(N\mathcal
{T})^{N^{3/2}}}{8^{N^{3/2}} \lceil N^{3/2} \rceil!}
\\[-2pt]
& \rightarrow & 0 \qquad\mbox{as } N \rightarrow\infty.
\end{eqnarray*}

Let $A_4'$ be the event that $\mathfrak{A}_t' = 0$ for some time $t
\in[T,T+\frac{1}{2}\mathcal{T})$. Notice that if $\mathfrak{A}_t' =
0$, then $\mathfrak{A}_s' = 0$ for $s \geq t$. Therefore, $A_4'$ is
the event that the label $\mathfrak{a}'$ is eliminated by time
$T+\frac{1}{2}\mathcal{T}$. By the given dynamics, $\mathfrak{A}_t'$
can only increase when individuals marked $\mathfrak{a}'$ replace
individuals marked $\mathfrak{b}'$ or $\mathfrak{c}'$ via resampling
events. At time $t$ the rate at which this happens is
%
%
\begin{equation}
\label{AUpper} \frac{1}{2} \cdot\frac{\mathfrak{A}_t'(N-\mathfrak
{A}_t')}{N} \leq
\mathfrak{A}_t'.
\end{equation}

We define the event $\mathcal{E}$ as
\[
\mathcal{E} = \bigl(A_1'\bigr)^C \cap
\bigl(A_2'\bigr)^C \cap A_3'
\cap\bigl\{T < \tfrac
{1}{2}\mathcal{T}\bigr\}.
\]
Selection will cause $\mathfrak{A}'$ to decrease. On the event
$(A_2')^C$ all of the individuals marked $\mathfrak{c}'$ will have
fitness at least $\frac{1}{32}W_0$ greater than any individual marked
$\mathfrak{a}$ until time $t_0$. Thus, on the event $(A_2')^C \cap\{T
< \frac{1}{2}t_0\}$, all of the individuals marked $\mathfrak{c}'$
will have fitness at least $\frac{1}{32}W_0$ greater than any
individual marked $\mathfrak{a}$ for all times $t \in[T,T+\frac
{1}{2}\mathcal{T})$. On the event $A_3'$ there are at least $N/5$
individuals marked $\mathfrak{c}$ for all times $t \in[T,T+\frac
{1}{2}\mathcal{T})$. Hence, on the event $\mathcal{E}$ individuals
marked $\mathfrak{a}'$ will become individuals marked $\mathfrak{c}'$
by a rate of at least
%
%
\begin{equation}
\label{ALower} \frac{\gamma\mathfrak{A}_t'\mathfrak{C}_t' W_0}{32N} \geq
\frac
{\gamma}{160} W_0
\mathfrak{A}_t'
\end{equation}
for all times $t \in[T,T+\frac{1}{2}\mathcal{T})$.

Let $\{U_n'\}$ be a biased random walk which goes up with probability
\[
p' = \frac{160}{160+\gamma W_0}
\]
and down with probability $1-p'$. Let $N$ be large enough so that $p' <
1/2$. Because the random walk is biased downward, the probability that
the random walk visits a state $j < U_0'$ is 1. Once the random walk is
in state $j$, it goes up 1 with probability $p'$ and will eventually
return to $j$ with probability~1. The random walk will go down 1 with
probability $1-p'$ and, from basic martingale arguments, the
probability that it never returns to $j$ again is $(1-2p')/(1-p')$.
Therefore, once $U'$ is in state $j$, the probability it never returns
to state $j$ is
\[
\frac{(1-2p')}{1-p'} \cdot\bigl(1-p'\bigr) = 1-2p'.
\]
Hence, the number of times $U'$ visits a state $j < U_0'$ has the
geometric distribution with mean $1/(1-2p')$. For more details see
\cite{DNA}, pages 194--196.\vadjust{\goodbreak}

By equations (\ref{AUpper}) and (\ref{ALower}) we see that on the
event $\mathcal{E}$, if $\mathfrak{A}'$ changes during the time
interval $[T,T+\frac{1}{2}\mathcal{T})$, it decreases with
probability higher than $p'$. The expected number of times that
$\mathfrak{A}'$ will visit state $j$ is therefore less than or equal
to $1/(1-2p')$ for any $j \in\{1,2,\ldots,N-1\}$. Also, the rate at
which $\mathfrak{A}_t'$ changes state is at least
\[
\frac{\gamma}{160}W_0\mathfrak{A}_t'
\]
for all times $t \in[T,T+\frac{1}{2}\mathcal{T})$ by equation (\ref
{ALower}). Let $\overline{A} = \{t \geq T\dvtx\mathfrak{A}_t' > 0\}$ and
let $\lambda$ be Lebesgue measure. Then
\[
E\bigl[\lambda(\overline{A})1_\mathcal{E}\bigr] \leq\frac
{160}{(1-2p')\gamma
W_0}
\sum_{j=1}^N \frac{1}{j} \sim
\frac{160\log N}{\gamma W_0}
\]
as $N \rightarrow\infty$.

Observe that
\begin{eqnarray*}
P\bigl(\mathcal{E} \cap\bigl(A_4'\bigr)^C
\bigr) & = & P \biggl(\mathcal{E} \cap\biggl\{ \lambda(\overline{A}) \geq
\frac{1}{2}\mathcal{T} \biggr\} \biggr)
\\
& = & P \biggl(\lambda(\overline{A})1_\mathcal{E} \geq\frac
{1}{2}
\mathcal{T} \biggr)
\\
& \leq &\frac{2 E[\lambda(\overline{A})1_\mathcal{E}]}{\mathcal{T}}
\qquad\mbox{by Markov's inequality}
\\
& \rightarrow &0 \qquad\mbox{as } N \rightarrow\infty.
\end{eqnarray*}
Therefore,
\[
P\bigl(\mathcal{E} \cap A_4'\bigr)-P \bigl(T<
\tfrac{1}{2}\mathcal{T} \bigr) \rightarrow0 \qquad\mbox{as } N \rightarrow
\infty.
\]

This allows us to do the following computation:
\begin{eqnarray*}
1 & = & \lim_{N \rightarrow\infty} \biggl(P \biggl(T < \frac
{1}{2}\mathcal{T}
\biggr)+P \biggl(T \geq\frac{1}{2}\mathcal{T} \biggr) \biggr)
\\
& = & \lim_{N \rightarrow\infty} \biggl(P\bigl(\mathcal{E} \cap A_4'
\bigr)+P \biggl(T \geq\frac{1}{2}\mathcal{T} \biggr) \biggr)
\\
& = & \lim_{N \rightarrow\infty} \biggl(P \biggl(\bigl(A_1'
\bigr)^C \cap\bigl(A_2'\bigr)^C
\cap A_3' \cap A_4' \cap
\biggl\{T < \frac{1}{2}\mathcal{T} \biggr\} \biggr)+P \biggl(T \geq
\frac{1}{2}\mathcal{T} \biggr) \biggr)
\\
& \leq &\lim_{N \rightarrow\infty} \biggl(P \biggl( \biggl\{B_T \leq
\frac{1}{2}\mathcal{T} \biggr\} \cap\biggl\{T < \frac{1}{2}\mathcal
{T} \biggr\} \biggr)+P \biggl(T \geq\frac{1}{2}\mathcal{T} \biggr)
\biggr)
\\
& = &\lim_{N \rightarrow\infty} P \biggl(B_T1_{\{T<
\mathcal{T}/2\}}\leq
\frac{1}{2}\mathcal{T} \biggr).
\end{eqnarray*}
\upqed
\end{pf}

Let $B = \inf\{t\dvtx X_t^- > X_0^+ - W_0/4\}$.
%
%
\begin{Prop} \label{BProp}
Suppose $W_0 \geq\mathcal{W}$ for all $N$. As $N$ tends to infinity,
\[
P(B > \mathcal{T}) \rightarrow0.
\]
\end{Prop}
\begin{pf}
First note that if $\mathfrak{B}_0+\mathfrak{C}_0 \geq N/4$ then,
because all of the individuals labeled $\mathfrak{b}$ or $\mathfrak
{c}$ at time 0 are also labeled $\mathfrak{c}'$, we have that
$\mathfrak{C}_0' \geq N/4$. The result then follows by Lemma \ref
{CLem} with $T = 0$. On the other hand, if $\mathfrak{B}_0+\mathfrak
{C}_0 < N/4$ then $\mathfrak{A}_0 \geq3N/4$.

Let $T = (\inf\{t\dvtx\mathfrak{A}_t < N/4\}) \wedge(\inf\{t\dvtx
\mathfrak{C}_t \geq N/4\})$. Let $A_5$ be the event that $\mathfrak
{A}_t \geq N/4$ for all times $t \in[0,\frac{1}{2}\mathcal{T})$. Let
$A_6$ be the event that $\mathfrak{C}_t < N/4$ for all times $t \in
[0,\frac{1}{2}\mathcal{T})$. Define $\zeta$ to be the infimum over
all times such that an individual labeled $\mathfrak{b}$ has fitness
in $(-\infty, X_0^+-\frac{5}{32}W_0)$, an individual labeled
$\mathfrak{c}$ has fitness in $(-\infty,X_0^+-\frac{3}{32}W_0)$ or
$\mathfrak{A}_t < N/4$. Note that $A_1^C \cap A_2^C \cap A_5 \subset\{
\zeta\geq\frac{1}{2}\mathcal{T}\}$.

On the event $\{\zeta\geq\frac{1}{2}\mathcal{T}\}$, the rate of
increase of $\mathfrak{C}_t$ due to selection is at least
%
%
\begin{equation}
\label{CLower} \frac{\gamma\mathfrak{A}_t \mathfrak{C}_t W_0}{32N} \geq
\frac
{1}{128}\gamma
\mathfrak{C}_t W_0
\end{equation}
for all $t \in[0,\frac{1}{2}\mathcal{T})$. On the other hand,
because $\mathfrak{C}_t$ can only decrease due to resampling,
$\mathfrak{C}_t$ will decrease no faster than
%
%
\begin{equation}
\label{CUpper} \frac{1}{2}\cdot\frac{\mathfrak{C}_t(N-\mathfrak
{C}_t)}{N} \leq
\mathfrak{C}_t.
\end{equation}

Let $\{U_n\}_{n=0}^\infty$ be a biased random walk with $U_0 = 1$
which goes up with probability
\[
p = \frac{\gamma W_0}{128+\gamma W_0}
\]
and down with probability $1-p$. Let $N$ be large enough so that $p >
1/2$. By similar reasoning as was used in the proof of Lemma \ref
{CLem}, the number of times $U_n$ visits a state $j \geq1$ has the
geometric distribution with mean $1/(2p-1)$. Also, by basic martingale
arguments, the probability that $U_n$ ever reaches state 0 is
\[
\frac{1-p}{p} = \frac{128}{\gamma W_0}.
\]

Note that $\mathfrak{C}_0 \geq U_0$ since the individual with the
highest fitness is initially labeled~$\mathfrak{c}$. On the event $\{
\zeta\geq\frac{1}{2}\mathcal{T}\}$, we see from equations (\ref
{CLower}) and (\ref{CUpper}) that if $\mathfrak{C}$ changes during
time $[0,\frac{1}{2}\mathcal{T})$, then it increases with a
probability of at least $p$. Therefore, the expected number of times
that $\mathfrak{C}$ visits state $j$ is less than or equal to
$1/(2p-1)$ and the probability the $\mathfrak{C}_t$ reaches state 0
for some time $t \in[0,\frac{1}{2}\mathcal{T})$ is less than
$128/(\gamma W_0)$. Let $A_7$ be the event that $\mathfrak{C}_t$
reaches state 0 for some time $t \in[0,\frac{1}{2}\mathcal{T})$.

By equation (\ref{CLower}), the rate at which $\mathfrak{C}$ changes
is at least
\[
\tfrac{1}{128}\gamma\mathfrak{C}_t W_0\vadjust{\goodbreak}
\]
for all times $t \in[0,\frac{1}{2}\mathcal{T})$ on the event $\{
\zeta> \frac{1}{2}\mathcal{T}\}$. Let $\overline{C} = \{t \in
[0,\frac{1}{2}\mathcal{T})\dvtx\mathfrak{C} < \frac{1}{4}N\}$ and let
$\lambda$ be Lebesgue measure. Then
\begin{eqnarray*}
E\bigl[\lambda(\overline{C})1_{\{\zeta\geq\mathcal{T}/2\}}\bigr] & = &
E\bigl[\lambda(
\overline{C})1_{\{\zeta\geq\mathcal{T}/2\}}1_{A_7}\bigr] + E\bigl
[\lambda(
\overline{C})1_{\{\zeta\geq\mathcal{T}/2\}}1_{A_7^C}\bigr]
\\
& \leq &\frac{1}{2}\mathcal{T} P(A_7) + \frac{128}{(2p-1)\gamma W_0}
\sum_{j=1}^{\lfloor N/4 \rfloor} \frac{1}{j}
\\
& \sim &\frac{128\log(N/4)}{\gamma W_0}.
\end{eqnarray*}

By Markov's inequality
\begin{eqnarray*}
P\bigl(A_1^C \cap A_2^C \cap
A_5 \cap A_6\bigr) & \leq & P\biggl(A_1^C
\cap A_2^C \cap A_5 \cap\biggl\{\lambda(
\overline{C}) \geq\frac{1}{2}\mathcal{T}\biggr\}\biggr)
\\
& \leq & P\biggl(\biggl\{\zeta\geq\frac{1}{2}\mathcal{T}\biggr\} \cap
\biggl
\{\lambda(\overline{C}) \geq\frac{1}{2}\mathcal{T}\biggr\}\biggr)
\\
& = & P\biggl(\lambda(\overline{C})1_{\{\zeta\geq\mathcal{T}/2\}} \geq
\frac{1}{2}
\mathcal{T}\biggr)
\\
& \leq &\frac{2E[\lambda(\overline{C})1_{\{\zeta\geq\mathcal{T}/2\}
}]}{\mathcal{T}}
\\
& \leq &\frac{256w^{1/4}\log(N/4)}{\mathcal{T}\gamma W_0} \qquad\mbox{for
}N\mbox{ large enough}
\\
& \rightarrow &0 \qquad\mbox{as } N \rightarrow\infty.
\end{eqnarray*}
Because $P(A_1^C \cap A_2^C) \rightarrow1$ we have $P(A_5^C \cup
A_6^C) \rightarrow1$ as $N \rightarrow\infty$.

Note that $A_5^C \cup A_6^C \subset\{T < \frac{1}{2}\mathcal{T}\}$.
Therefore, $P(T < \frac{1}{2}\mathcal{T}) \rightarrow1$ as $N
\rightarrow\infty$. Let $E_2 = (A_1')^C \cap(A_2')^C \cap\{T <
\frac{1}{2}\mathcal{T}\}$. Then $P(E_2) \rightarrow1$ as $N
\rightarrow\infty$. To show $P(B \leq\mathcal{T}) \rightarrow1$ we
can show $P(\{B \leq\mathcal{T}\} \cap E_2) \rightarrow1$. At time
$T$, at least $\frac{1}{4}N$ individuals will be labeled either
$\mathfrak{b}$ or $\mathfrak{c}$. According to the labeling, all of
these individuals are labeled $\mathfrak{c}'$ so that at time $T$ we
have $\mathfrak{C}_T \geq\frac{1}{4}N$. By Lemma~\ref{CLem} we have
\[
P \bigl(B_T 1_{\{T < \mathcal{T}/2\}} \leq\tfrac
{1}{2}\mathcal{T}
\bigr) \rightarrow1 \qquad\mbox{as } N \rightarrow\infty.
\]
Note that
\[
\bigl\{B_T 1_{\{T < \mathcal{T}/2\}} \leq\tfrac
{1}{2}\mathcal{T}
\bigr\} = \bigl\{B_T \leq\tfrac{1}{2}\mathcal{T} \bigr\}
\cup\bigl\{T \geq\tfrac{1}{2}\mathcal{T} \bigr\}.
\]
Because $E_2 \subset\{T < \frac{1}{2}\mathcal{T}\}$ we have
\[
\bigl\{B_T 1_{\{T < \mathcal{T}/2\}} \leq\tfrac
{1}{2}\mathcal{T}
\bigr\} \cap E_2 = \bigl\{B_T \leq\tfrac
{1}{2}
\mathcal{T} \bigr\} \cap E_2.
\]
It then follows that
\[
P \bigl( \bigl\{B_T \leq\tfrac{1}{2}\mathcal{T} \bigr\}
\cap E_2 \bigr) \rightarrow1 \qquad\mbox{as } N \rightarrow\infty.
\]
However,
\[
\bigl\{B_T \leq\tfrac{1}{2}\mathcal{T} \bigr\} \cap
E_2 \subset\bigl\{B_T \leq\tfrac{1}{2}
\mathcal{T} \bigr\} \cap\bigl\{T < \tfrac{1}{2}\mathcal{T} \bigr\}
\subset\{B \leq\mathcal{T}\},
\]
which gives the conclusion.
\end{pf}

Let $V_t^1 = \{i\dvtx X_t^i > X_0^+ + W_0/4\}$ and $V_t^2 = \{i\dvtx
X_t^i <
X_0^- - W_0/4\}$. Let $F = \inf\{t\dvtx V_t^1 \cup V_t^2 \neq\varnothing
\}
$. We now want to bound the time it takes for the width to increase.
%
%
\begin{Prop} \label{FProp}
Suppose $W_0 \geq\mathcal{W}$ for all $N$. Then
\[
\lim_{N \rightarrow\infty} P(F > \mathcal{T}) = 1.
\]
\end{Prop}
\begin{pf}
By Proposition~\ref{UpBound} with $l = W_0/4$ and $t = \mathcal{T}$
we have
\begin{eqnarray*}
P\bigl(\inf\bigl\{s\dvtx V_s^1 \neq\varnothing\bigr\} < t
\bigr) & = & P \Bigl(\sup_{0 \leq s
\leq t} D_s \geq l \Bigr)
\\
& \leq &\frac{2N(t\mu)^le^{(\gamma(W_0+2l)+\mu+1)t}}{(l-1)!}
\\
& \rightarrow & 0 \qquad\mbox{as } N \rightarrow\infty.
\end{eqnarray*}

By Proposition~\ref{BackSpeed} with $l = W_0/4$ and $t = \mathcal{T}$
we have
\begin{eqnarray*}
P\bigl(\inf\bigl\{s\dvtx V_s^2 \neq\varnothing\bigr\} < t
\bigr) & = & P \Bigl(\sup_{0 \leq s
\leq t} \bigl(X_0^- -
X_s^-\bigr) \geq l \Bigr)
\\
& \leq &\frac{N(t\mu)^l e^t}{l!}
\\
& \rightarrow & 0 \qquad\mbox{as } N \rightarrow\infty.
\end{eqnarray*}
\upqed
\end{pf}

Recall that\vspace*{1pt} $Y_i = \sup_{s_i \leq s \leq t_{i+1}}D_{s} -
D_{s_i}$ and that $\{\mathcal{F}_t\}_{t \geq0}$ is the natural
filtration associated with $X$. Note that if $W_0 < 2\mathcal{W}$, then
for all $n \geq1$ the width satisfies $W_{s_n} = \lceil2\mathcal{W}
\rceil$.
\begin{pf*}{Proof of Proposition~\ref{TheoProp2}}
We consider a sequence of initial configurations $X_0$ depending on $N$
such that $W_0 = \lceil2\mathcal{W} \rceil$ for all $N$. Because
$W_0 \geq2\mathcal{W}$ we have $s_1 = 0$ and $Y_1 = \sup_{0 \leq s
\leq t_2} D_s-D_0$. We will show that for $N$ large enough, $E[Y_1] <
5\mathcal{W}$. The result then follows because $X$ is a strong Markov process.

We make the following definitions:
\begin{eqnarray*}
V_t^1(s) &=& \bigl\{i\dvtx X_t^i >
X_s^+ + W_s/4\bigr\} \qquad\mbox{for } t \geq s \geq0,
\\
V_t^2(s) &=& \bigl\{i\dvtx X_t^i <
X_s^- - W_s/4\bigr\} \qquad\mbox{for } t \geq s \geq0,
\\
F_0 &=& B_0 = r_0 = 0,
\\
F_n &=& \inf\bigl\{t \geq r_{n-1}\dvtx V_t^1(r_{n-1})
\cup V_t^2(r_{n-1}) \neq\varnothing\bigr\}
\qquad\mbox{for }n \geq1,
\\
B_n &=& \inf\bigl\{t \geq r_{n-1}\dvtx X_t^- >
X_{r_{n-1}}^+-W_{r_{n-1}}/4\bigr\} \qquad\mbox{for }n \geq1,
\\
r_n &=& F_n \wedge B_n \qquad\mbox{for }n \geq1,
\\
n_* &=& \inf\{n \geq1\dvtx W_{r_n} < \mathcal{W}\}.
\end{eqnarray*}
Note that $r_1$ is the first time that the event $F \cup B$ occurs and
that, conceptually, $r_n$~acts like the first time that $F \cup B$
occurs when the process is started at time $r_{n-1}$ for $n \geq2$.
The random variables $F_n$ and $B_n$ play the roles of the events $F$
and $B$ when the processes are started at time $r_{n-1}$.

On the event $n-1 < n_*$, by Proposition~\ref{BProp} and the strong
Markov property of $X$, we have $P(B_n \leq r_{n-1}+\mathcal
{T}|\mathcal{F}_{r_{n-1}}) \rightarrow1$ uniformly on a set of
probability 1 as $N \rightarrow\infty$. Likewise, on the event $n-1 <
n_*$, by Proposition~\ref{FProp} and the strong Markov property, we
have $P(F_n > r_{n-1}+\mathcal{T}|\mathcal{F}_{r_{n-1}}) \rightarrow
1$ uniformly on a set of probability 1 as $N \rightarrow\infty$.
Therefore, on the event $n-1 < n_*$, we have $P(B_n < F_n|\mathcal
{F}_{r_{n-1}}) \rightarrow1$ uniformly on a set of probability 1.

Because the bounds in Propositions~\ref{BProp} and~\ref{FProp} do not
depend on $n$ we can choose a sequence $p = p_N$ such that $p
\rightarrow1$ as $N \rightarrow\infty$ and almost surely
\[
p 1_{\{n-1 < n_*\}} \leq P(B_n < F_n|
\mathcal{F}_{r_{n-1}})1_{\{n-1 <
n_*\}}
\]
for all $n \geq0$. Let $\{S_n\}_{n=0}^\infty$ be a random walk
starting at 1 which goes down 1 with probability $p$ and up 1 with
probability $1-p$ until it reaches 0. Once $S$ reaches 0 it is fixed.
For $n < n_*$ we couple $S$ with $X$ so that $2^{S_n-1}W_0 \geq
W_{r_n}$. The coupling is defined as follows:
\begin{itemize}
\item Each step of the process $S$ corresponds to a time $r_n$.
\item On the event $\{F_n < B_n\}$ we have $S_n - S_{n-1} = 1$.
\item On the event $\{B_n \leq F_n\}$ we have $S_n - S_{n-1} = -1$ with
probability $p/P(B_n \leq F_n)$ and we have $S_n - S_{n-1} = 1$ with
probability $1-p/P(B_n \leq F_n)$.
\end{itemize}

We will\vspace*{1pt} show that this coupling is well defined and gives
the necessary bound. Initially, $S_0 = 1$ and $2^{S_0-1}W_0 = W_0$. On
the event that $B_n \leq F_n$, we have $W_{r_n} <
\frac{1}{2}W_{r_{n-1}}$ and $\sup _{r_{n-1} \leq t \leq r_n}
D_t-D_{r_{n-1}} \leq\frac{1}{4} W_{r_{n-1}}$. On the event that $F_n <
B_n$, we have $W_{r_n} < 2W_{r_{n-1}}$ and $\sup_{r_{n-1} \leq t \leq
r_n}D_t-D_{r_{n-1}} \leq \frac{1}{4}W_{r_{n-1}}+1$. Therefore, if
$2^{S_{n-1}-1}W_0 \geq W_{r_{n-1}}$, then $2^{S_n-1}W_0 \geq W_{r_n}$
by the coupling. It follows that $2^{S_n-1}W_0 \geq\sup_{r_{n-1} \leq t
\leq r_n} D_t-D_{r_{n-1}}$ as well. By induction,\vspace*{1pt}
$2^{S_n-1}W_0 \geq W_{r_n}$ for all $n < n_* \wedge\inf\{m\dvtx S_m =
0\}$. If $n = \inf\{m\dvtx S_m=0\}$, then $W_{r_n} \leq\mathcal{W}$.
Therefore, $n_* \leq\inf\{m\dvtx S_m=0\} $ and the induction holds for
all $n < n_*$.

We define a function $d$ on $(\{0\} \cup\N)^\infty$ such that if $x
= (x_0, x_1, \ldots)$ then
\[
d(x) = \sum_{i=0}^\infty1_{\{x_i > 0\}}
2^{x_i-1} W_0.
\]
Consider $S = (S_0, S_1, \ldots)$ as a random element in $(\{0\} \cup
\N)^\infty$. Then
\[
d\bigl((S_0, S_1, \ldots, S_n, 0, 0, \ldots)
\bigr) \geq\sum_{i=1}^n \Bigl(
\sup_{r_{i-1} \leq t \leq r_i} D_t-D_{r_{i-1}} \Bigr) \geq
\sup_{0 \leq t
\leq r_n} D_t
\]
for all $n$ such that $n-1 < n_*$. By definition, $n_*$ is the first
$n$ such that $W_{r_n} < \mathcal{W}$. Hence, $d(S) \geq
Y_1$.

For any $n \geq0$ we have
\[
P(S_{2n+1} = 0) = \pmatrix{2n+1
\cr
n}(1-p)^n
p^{n+1} \leq4^n (1-p)^n p^{n+1}.
\]
If $S_{2n+1} = 0$ then
\[
d(S) \leq\Biggl(2 + 2\sum_{i=1}^n
2^{i-1} \Biggr) W_0 = 2^{n+1} W_0,
\]
which is obtained by taking $n$ steps up followed by $n+1$ steps down.

Therefore,
\[
E[Y_1] \leq E\bigl[d(S)\bigr] \leq\sum_{n=0}^\infty
\bigl[4(1-p)\bigr]^n p^{n+1} 2^{n+1}W_0
= \frac{2p W_0}{1-8(1-p)p} \sim4\mathcal{W},
\]
because $W_0 = \lceil2\mathcal{W} \rceil$ and $p \rightarrow1$ as
$N \rightarrow\infty$. This shows that for $N$ large enough we have
$E[Y_1] < 5\mathcal{W}$, which gives the conclusion.
\end{pf*}

Let $l = \lfloor\mathcal{W}/2 \rfloor$. We make the following
definitions for the rest of the section:
\begin{eqnarray*}
K_1 &=& \frac{2N(\mathcal{T}\mu)^l e^{(\gamma(W_0+2l)+\mu
+1)\mathcal{T}}}{(l-1)!},
\\[-2pt]
K_2 &=& \frac{N(\mathcal{T}\mu)^l e^{\mathcal{T}}}{l!},
\\[-2pt]
p &=& 1-K_1-K_2.
\end{eqnarray*}

%
\begin{Lemma} \label{WidthSpeed}
Suppose $W_0 \leq\mathcal{W}$ for all $N$. Then
\[
P \Bigl(\sup_{0 \leq s \leq\mathcal{T}} W_s \leq2\mathcal{W} \Bigr)
\geq1-K_1-K_2.
\]
\end{Lemma}
\begin{pf}
By Proposition~\ref{UpBound} we have
\[
P \Bigl(\sup_{0 \leq s \leq\mathcal{T}} D_s \geq l \Bigr) \leq
K_1.
\]
By Proposition~\ref{BackSpeed} we have
\[
P \Bigl(\sup_{0 \leq s \leq\mathcal{T}} \bigl(X_0^- - X_s^-\bigr)
\geq l \Bigr) \leq K_2.
\]
On the event that $\sup_{0 \leq s \leq t} D_s \leq\mathcal{W}/2$ and
$\sup_{0 \leq s \leq t}X_0^--X_s^- \leq\mathcal{W}/2$, we have $\sup_{0
\leq s \leq t} W_t \leq2\mathcal{W}$. This gives the result.\vadjust{\goodbreak}
\end{pf}
\begin{pf*}{Proof of Proposition~\ref{TheoProp3}}
Notice that
\[
\{N_s \geq i\} = \{s_i \leq s\} \subset\Biggl\{\sum
_{j=1}^i (s_j-t_j)
\leq s \Biggr\}.
\]
Therefore,
\[
P(N_s \geq i) \leq P \Biggl(\sum_{j=1}^i
(s_j-t_j) \leq s \Biggr).
\]
Applying Lemma~\ref{WidthSpeed} and the strong Markov property of $X$
we have
\[
1-K_1-K_2 \leq P (s_j-t_j \geq
\mathcal{T}|\mathcal{F}_{t_j} )
\]
for all $j$. Taking expectations of both sides yields
\[
1-K_1-K_2 \leq P(s_j-t_j \geq
\mathcal{T})
\]
for all $j$, so
\[
1-K_1-K_2 \leq\inf_j P(s_j-t_j
\geq\mathcal{T}).
\]

Note that $p \rightarrow1$ as $N \rightarrow\infty$. Define an
i.i.d. sequence $\{V_i\}_{i=1}^\infty$ of random variables with
distribution $P(V_i = 0) = 1-p$ and $P(V_i = \mathcal{T}) = p$. Then
\[
P \Biggl(\sum_{j=1}^i
(s_j-t_j) \leq s \Biggr) \leq P \Biggl(\sum
_{j=1}^i V_i \leq s \Biggr).
\]
This will allow us to define a new process $N_s'$ such that $N_s' = i$ if
\[
\sum_{j=1}^i V_i \leq s <
\sum_{j=1}^{i+1} V_i.
\]
Note that $P(N_s' = 0) = p$ for $s \in[0,\mathcal{T})$ and that
$P(N_s' \geq k) \geq P(N_s \geq k)$ for all~$k$. Therefore, it is
enough to bound $E[N_s']/s$.

Let $V_0 = 0$. Jumps of the process $N_s'$ only occur at points
$k\mathcal{T}$ where $k$ is a positive integer. On the time interval
$[0,\mathcal{T})$ the process $N_s'$ is constant and has value $\max\{
i \geq0\dvtx V_i = 0\}$. Therefore, $N_s'$ has the shifted geometric
distribution for $s \in[0,\mathcal{T})$ with mean $(1-p)/p$. We can
now make use of the fact that $N_s'$ is a Markov process. If we
consider values at $k\mathcal{T}$ for $k \geq0$, we have for $s \in
[(k-1)\mathcal{T}, k\mathcal{T})$ that $E[N_s'] = k(1-p)/p$. For $k
\geq2$ we then have
\[
\frac{1}{s}E\bigl[N_s'\bigr] =
\frac{k(1-p)}{sp} \leq\frac{k(1-p)}{(k-1)p
\mathcal{T}}.
\]
This gives us
\[
\frac{\mathcal{T}}{s}E\bigl[N_s'\bigr] \leq
\frac{k(1-p)}{(k-1)p} \rightarrow0 \qquad\mbox{as }N \rightarrow\infty.
\]
On the time interval $[0,\mathcal{T})$ we have
\[
\frac{\mathcal{T}}{s}E\bigl[N_s'\bigr] \leq
\frac{(1-p)}{p} \rightarrow0 \qquad\mbox{as } N \rightarrow\infty.
\]
\upqed
\end{pf*}

\begin{center}
\begin{tabular*}{\tablewidth}{@{\extracolsep{\fill}}l p{0.85\textwidth}@{}}
\multicolumn{2}{@{}c@{}}{NOTATION}\\[12pt]
$N$ & The size of the population \\
$\mu$ & The rate at which individuals accumulate mutations \\
$q$ & The probability that a mutation is beneficial \\
$\gamma$ & The selection coefficient \\
$X^i$ & The stochastic process in $\Z$ that represents the fitness of
the $i$th individual \\
$X$ & The stochastic process in $\Z^N$ that represents the fitnesses
of the individuals \\[2.5pt]
$\overline{X}$ & $= \frac{1}{N}\sum_{i=1}^N X^i$ \\[2.5pt]
$X_t^+$ & $= \max\{X_t^i\dvtx1 \leq i \leq N\}$ \\[2.5pt]
$X_t^-$ & $= \min\{X_t^i\dvtx1 \leq i \leq N\}$ \\[2.5pt]
$W_t$ & $= X_t^+-X_t^-$ \\[2.5pt]
$D_t$ & $= X_t^+ - X_0^+$ \\[2.5pt]
$w$ & is any positive, increasing function satisfying $\lim_{N
\rightarrow\infty} w(N) = \infty$ \\[2.5pt]
& and $\lim_{N \rightarrow\infty} w(N)/\log\log N = 0$\\[2.5pt]
$\mathcal{W}$ & $= \lfloor w\log N/\log\log N \rfloor$ \\[2.5pt]
$\mathcal{T}$ & $= w^{-1/2}\log\log N$ \\[2.5pt]
$t_1$ & $= 0$ \\[2.5pt]
$s_n$ & $= \inf\{t \geq t_n\dvtx W_t \geq2\mathcal{W}\} \mbox{ for } n
\geq1$\\[2.5pt]
$t_n$ & $= \inf\{t \geq s_{n-1}\dvtx W_t < \mathcal{W}\} \mbox{ for } n
\geq2$ \\[2.5pt]
$Y_i$ & $= \sup_{s_i \leq t \leq t_{i+1}}D_t - D_{s_i} \mbox{ for }i
\geq1$ \\[2.5pt]
$N_t$ & $= \max\{i\dvtx s_i \leq t\} \mbox{ for }t \geq0$ \\[2.5pt]
$Z_t^{k,\uparrow}$ & A multi-type Yule process in which there are
initially $N$ particles of type $k$.
Particles increase from type $i$ to type $i+1$ at rate $\mu$ and
particles of type $i$
branch at rate $\gamma i + 1$ \\
$\overline{M}{}^{k,\uparrow}_t$ & The maximum type of any particle in
$Z_t^{k,\uparrow}$ \\[2.5pt]
$M_t^{k,\uparrow}$ & $\overline{M}{}^{k,\uparrow}_t-k$\\[2.5pt]
$X_t'$ & $X_0^++\mathcal{M}_t^0$ if $t \in[0,\mathcal{T}]$ and
$X_{i\mathcal{T}}'+\mathcal{M}_{t-i\mathcal{T}}^i$ if $t \in
(i\mathcal{T}, (i+1)\mathcal{T}]$ for any\\[2pt]
$\{\mathcal{Z}_t^n\}_{n=0}^\infty$ & An i.i.d.\vspace*{2pt} sequence of stochastic
processes each having the same distribution
as $Z^{\mathcal{W},\uparrow}$ \\[2.5pt]
$\overline{\mathcal{M}}{}^n_t$ & The maximum type of any particle in
$Z^{\mathcal{W},\uparrow}$
\end{tabular*}
\end{center}

\begin{center}
\begin{tabular*}{\tablewidth}{@{\extracolsep{\fill}}l p{0.85\textwidth}@{}}
$\mathcal{M}_t^n$ & $=\overline{\mathcal{M}}{}^n_t - \mathcal{W}$\\[2.5pt]
& integer $i \geq1$ \\[2.5pt]
$D_t'$ & $X_t' - X_0^+$ \\[2pt]
$\mathcal{F}$ & $=\{\mathcal{F}_t\}_{t \geq0}$ is the natural
filtration associated with $X$ under the initial condition
$X_0^i = 0$ for $1 \leq i \leq N$ \\[2pt]
$Z_t^C$ & A multi-type Yule process in which there are initially $N$
particles of type 0.
Particles increase from type $i$ to type $i+1$ at rate $\mu$ and
branch at rate~$C$ \\[2pt]
$M_t^C$ & The maximum type of any particle in $Z_t^C$ \\[2pt]
$S_t$ & $=\sup_{0 \leq s \leq t}(X_0^- - X_s^-)$\\[2pt]
$A_t^{x,l}$ & The event that an individual in $\overline{S}{}^x_s$ has
fitness in $(-\infty, x-l]$
for some time $s \in[0,t]$\\[2pt]
$A_1$ & The event that there is an
individual labeled $\mathfrak{b}$ with fitness in $(-\infty,
X_0^+-\frac{5}{32}W_0)$
for some time $t \in[0,\mathcal{T}]$\\[2pt]
$A_2$ & The event that there is an
individual labeled $\mathfrak{c}$ with fitness in $(-\infty,X_0^+ -
\frac{3}{32}W_0)$
for some time $t \in[0,\mathcal{T}]$\\[2pt]
$A_1'$ & The event that there is an
individual labeled $\mathfrak{b}'$ with fitness in $(-\infty,
X_0^+-\frac{7}{32}W_0)$
for some time $t \in[0,\mathcal{T}]$\\[2pt]
$A_2'$ & The event that there is an
individual labeled $\mathfrak{c}'$ with fitness in $(-\infty,X_0^+ -
\frac{5}{32}W_0)$
for some time $t \in[0,\mathcal{T}]$\\[2pt]
$B$ & $= \inf\{t\dvtx X_t^- > X_0^+ -
W_0/4\}$\\[2pt]
$V_t^1$ & $=\{i\dvtx X_t^i > X_0^++W_0/4\}
$\\[2pt]
$V_t^2$ & $=\{i\dvtx X_t^i < X_0^--W_0/4\}
$\\[2pt]
$F$ & $=\inf\{t\dvtx V_t^1 \cup V_t^2 \neq
\varnothing\}$
\end{tabular*}
\end{center}

\begin{appendix}\label{app}
\section*{Appendix}
%
%
\begin{Lemma} \label{TRTBoundLem}
Let $x \geq0$. The tail of the exponential series satisfies
\[
\sum_{i=k}^\infty\frac{x^i}{i!} \leq
\frac{x^k e^x}{k!}.\vspace*{-2pt}
\]
\end{Lemma}
\begin{pf}
By Taylor's remainder theorem we know that there exists a $\xi\in
[0,x]$ such that
\[
e^x = \sum_{i=1}^{k-1}
\frac{x^i}{i!} + \frac{x^k e^{\xi}}{k!}.
\]
Using the series expansion of $e^x$ we have
\[
\sum_{i=k}^\infty\frac{x^i}{i!} =
\frac{x^k e^{\xi}}{k!} \leq\frac{x^k e^x}{k!}.\hspace*{125pt} \qed\hspace*{-125pt}
\]
\noqed\end{pf}

Recall that $M_t^C$ is the maximum type of any particle in the
branching process~$Z_t^C$.
%
%
\begin{Lemma} \label{YuleLem}
For any population size $N$, time $t \geq0$ and natural number~$l$,
\[
P\bigl(M_t^C \geq l\bigr) \leq\frac{N (t\mu)^l e^{Ct}}{l!}.\vspace*{-2pt}
\]
\end{Lemma}
\begin{pf}
Consider a Yule process $Z$ which is the same as $Z^C$ except there is
only one particle at time 0. It is well known that the number of
particles in $Z_t$ has mean $e^{Ct}$. Let $M_t'$ be the maximum type of
any particle at time $t$. When there are $k$ particles in the
population, we let $B_1, \ldots, B_k$ denote the types of the
particles, where the numbering is independent of the mutations. For any
$l \geq0$,
\begin{eqnarray*}
P\bigl(M_t' \geq l\bigr) & = & \sum
_{k=1}^\infty P\bigl(M_t'
\geq l|Z_t = k\bigr)P(Z_t = k)
\\[-2pt]
& = & \sum_{k=1}^\infty P\bigl(
\{B_1 \geq l\} \cup\cdots\cup\{B_k \geq l\}
|Z_t = k\bigr)P(Z_t = k)
\\[-2pt]
& \leq &\sum_{k=1}^\infty k P(B_1
\geq l)P(Z_t = k)
\\[-2pt]
& = & E[Z_t]P(B_1 \geq l)
\\[-2pt]
& = & e^{Ct} \sum_{i=l}^\infty
\frac{(t\mu)^i}{i!}e^{-\mu t}.
\end{eqnarray*}
By Lemma~\ref{TRTBoundLem} it follows that
\[
P\bigl(M_t' \geq l\bigr) \leq\frac{(t\mu)^l e^{Ct}}{l!}.
\]

Now consider\vspace*{1pt} $Z^C$. At time 0 label the particles $1,2,\ldots,N$ and
let $M_{i,t}'$ be the maximum type of any particle among the progeny of
particle $i$ at time~$t$. Then
\begin{eqnarray*}
P\bigl(M_t^C \geq l\bigr) & = & P\bigl(\bigl
\{M_{1,t}' \geq l\bigr\} \cup\cdots\cup\bigl
\{M_{N,t}' \geq l\bigr\}\bigr)
\\[-2pt]
& \leq & NP\bigl(M_{1,t}' \geq l\bigr)
\\[-2pt]
& \leq &\frac{N (t\mu)^l e^{Ct}}{l!}.
\end{eqnarray*}
\upqed\end{pf}

Recall that $M_t^{k,\uparrow} = \overline{M}{}^{k,\uparrow}_t-k$ where
$\overline{M}{}^{k,\uparrow}_t$ is the maximum type of any individual
in the branching process $Z_t^{k,\uparrow}$.
%
%
\begin{Lemma} \label{ZUpBound}
For any time $t \geq0$ and any integers $k \geq0$ and $l \geq0$ we~have
\[
P\bigl(M_t^{k,\uparrow} > l\bigr) \leq\frac{N (t\mu)^l e^{(\gamma
(k+l)+1)t}}{l!}.
\]
\end{Lemma}
\begin{pf}
While all of the particles in $Z_t^{k,\uparrow}$ have type less than
$k+l$, they branch at a rate which is less than or equal to $\gamma
(k+l)+1$. Because of this, $P(M_t^{k,\uparrow} > l) \leq P(M_t^{\gamma
(k+l)+1} > l)$. By Lemma~\ref{YuleLem} we have
\[
P\bigl(M_t^{\gamma(k+l)+1} > l\bigr) \leq\frac{N (t\mu)^l e^{(\gamma
(k+l)+1)t}}{l!}.
\]
\upqed\end{pf}
\end{appendix}

\section*{Acknowledgments}

I would like to thank Jason Schweinsberg for suggesting the problem,
patiently helping me work through various parts of the proof and for
helping me revise the first drafts of the paper. I would also like to
thank the referee for helpful comments that led to an improved upper
bound.


%

\printaddresses

\end{document}